\documentclass[11pt,a4paper]{article}
\usepackage[a4paper, total={6in, 8in}]{geometry}

\usepackage[utf8]{inputenc}
\usepackage{hyperref}
\usepackage{tabularx}
\usepackage{enumitem} 
\newlist{condenum}{enumerate}{1} 
\setlist[condenum]{label=\bfseries P\arabic*., 
                   ref=\arabic*, wide}

\usepackage[cmex10]{amsmath}
\usepackage{%
	amsfonts,%
	amsmath,%
	amssymb,%
	amsthm,%
	bbm,%
	cite,%
	float,%
	etoolbox,%
	mathrsfs,%
	mathtools,%
	multirow,%
	pgf,%
	pgfplots,%
	subfigure,%
	tikz,%
}
\usepackage[normalem]{ulem}
\usepackage{cleveref}
\pgfdeclareplotmark{mystar}{
    \node[star,star point ratio=2.25,minimum size=6pt,
          inner sep=0pt,draw=black,solid,fill=red] {};
}
\pgfdeclareplotmark{mystar2}{
    \node[star,star point ratio=2.25,minimum size=6pt,
          inner sep=0pt,draw=black,solid,fill=blue] {};
}

\newlist{temp}{enumerate}{1} 
\setlist[temp]{label=(\roman*), 
                   ref=\roman*, wide}
\usepgflibrary{shapes}
\usetikzlibrary{%
	arrows,%
	patterns,%
	positioning,%
	shapes.callouts,%
	shapes.symbols,%
}
\usepackage[ruled,vlined]{algorithm2e}
\usepackage{bbm}
\usepackage{breqn}

\newcommand\numeq[1]%
  {\stackrel{\scriptscriptstyle(\mkern-1.5mu#1\mkern-1.5mu)}{=}}

\theoremstyle{plain}
\newtheorem{thm}{Theorem}
\newtheorem{lem}[thm]{Lemma}
\newtheorem{proposition}[thm]{Proposition}
\newtheorem{cor}[thm]{Corollary}

\theoremstyle{definition}
\newtheorem{defn}[thm]{Definition}

\newtheoremstyle{case}{}{}{}{}{}{:}{ }{}
\theoremstyle{case}

\theoremstyle{remark}

\newcommand{\eq}[1]{\begin{align*}#1\end{align*}}
\newcommand{\EQ}[1]{\begin{equation*}#1\end{equation*}}
\newcommand{\eqn}[1]{\begin{align}#1\end{align}}
\newcommand{\EQN}[1]{\begin{equation}#1\end{equation}}

\newcommand{\norm}[1]{\left\lVert#1\right\rVert}

\DeclarePairedDelimiter\abs{\lvert}{\rvert}%

\newcommand{\Nats}{\mathbb{N}}
\newcommand{\E}{\mathbb{E}}
\newcommand{\Z}{\mathbb{Z}}

\newcommand{\R}{\mathbb{R}}

\renewcommand{\ge}{\geqslant}
\renewcommand{\le}{\leqslant}
\newcommand{\floor}[1]{\lfloor #1 \rfloor}
\newcommand{\mb}[1]{\mathbb{#1}}
\newcommand{\mc}[1]{\mathcal{#1}}
\newcommand{\mf}[1]{\mathbf{#1}}

\newcommand{\brac}[1]{\left(#1\right)}
\newcommand{\cbrac}[1]{\left\{#1\right\}}
\newcommand{\sbrac}[1]{\left[#1\right]}
\newcommand{\indic}[1]{\mathbbm{1}{\brac{#1}}}
\newcommand{\expect}[1]{\mathbb{E}\sbrac{#1}}
\newcommand{\define}{\triangleq}

\title{Asymptotic Optimality of Speed-Aware JSQ for Heterogeneous Service Systems}
\author{Sanidhay Bhambay and Arpan Mukhopadhyay\\
University of Warwick}
\date{}

\begin{document}
\maketitle
\begin{abstract}
    The Join-the-Shortest-Queue (JSQ) load-balancing scheme is known to minimise the average delay of jobs in homogeneous systems consisting of identical servers. However, it performs poorly in heterogeneous systems where servers have different processing rates. Finding a delay optimal scheme remains an open problem for heterogeneous systems. In this paper, we consider a speed-aware version of the JSQ scheme for heterogeneous systems and show that it achieves delay optimality in the fluid limit. One of the key issues in establishing this optimality result for heterogeneous systems is to show that the sequence of steady-state distributions indexed by the system size is tight in an appropriately defined space. The usual technique for showing tightness by coupling with a suitably defined dominant system does not work for heterogeneous systems. To prove tightness, we devise a new technique that uses the drift of exponential Lyapunov functions. Using the non-negativity of the drift, we show that the stationary queue length distribution has an exponentially decaying tail - a fact we use to prove tightness. Another technical difficulty arises due to the complexity of the underlying state-space and the separation of two time-scales in the fluid limit.  Due to these factors, the fluid-limit turns out to be a function of the invariant distribution of a multi-dimensional Markov chain which is hard to characterise. 
    By using some properties of this invariant distribution and using the monotonicity of the system, we show that the fluid limit is has a unique and globally attractive fixed point.
\end{abstract}

\section{Introduction}

The average response time of user requests is a key performance measure in modern large-scale service systems such as web server farms and cloud data centers. How incoming user requests or {\em jobs} are assigned to {\em servers}  can significantly impact the performance of such systems. The canonical model for studying the effect of job assignment schemes on the mean response time of jobs consists of $N$ servers each having its own queue, a stream of incoming jobs arriving at rate $N\lambda$, and a job dispatcher. The job dispatcher assigns every incoming job to a server in the system according to a specific job assignment scheme.
For homogeneous systems consisting of identical servers, a natural job assignment scheme
to consider is the Join-the-Shortest-Queue (JSQ) scheme in which each incoming job is assigned to the server with the minimum queue length. For finite $N$, this scheme is known to achieve the minimum average response time of jobs in a variety of settings~\cite{Winston1977optimality,Weber1978,johri1989optimality,hordijk1990optimality}.
Furthermore, in the large system limit ($N \to \infty$) it has been shown that the average queuing delay of jobs (average time a job spends in the queue before its processing starts) under JSQ approaches to zero. 

Variants of the JSQ scheme such as SQ($d$) scheme and the Join-the-Idle-Queue (JIQ) scheme have also been extensively studied in the literature~\cite{Lu2011,Gamarnik2016delay,Mukherjee2018}. In the SQ($d$) scheme, each job is assigned to the shortest of $d$ randomly sampled queues; in the JIQ scheme, each job is assigned to an idle queue if one is available and is otherwise sent to a randomly sampled queue. These schemes require less messaging overhead than JSQ but achieve the same asymptotic performance as JSQ in the large system limit. 

\begin{figure}[h!]
  \centering
  \includegraphics[width=8cm]{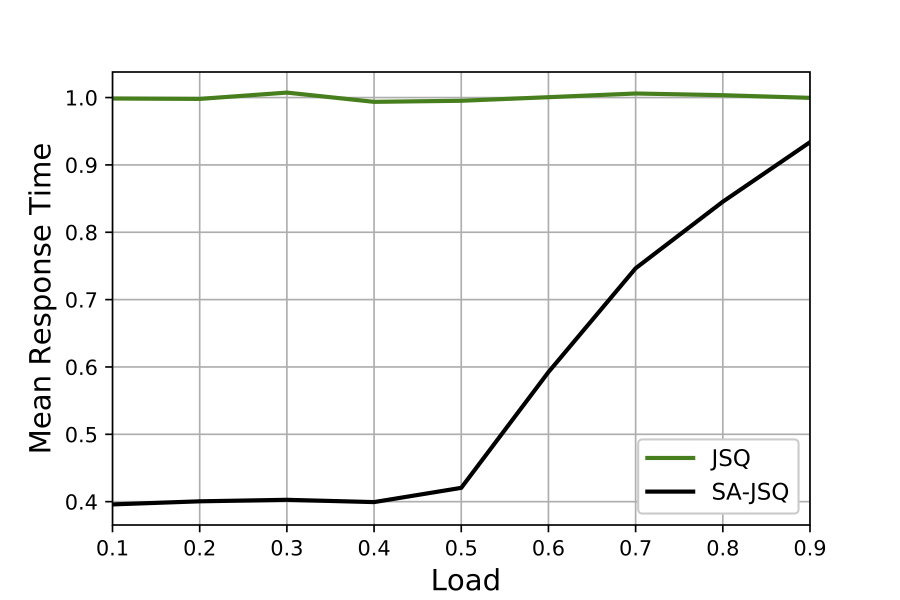}
\caption{Comparison of mean response time of jobs under JSQ and SA-JSQ schemes. 
Both schemes are applied to a heterogeneous system consisting of two types of servers:
a fraction $\gamma_i$ of servers have rate $\mu_i$. We choose
$\gamma_1=1-\gamma_2=1/5$, and $\mu_1=4\mu_2=20/8$.}
\label{fig:1}
\end{figure}

The (asymptotic) optimality of JSQ and its variants crucially relies on the assumption
of homogeneity of server speeds.  However, this assumption is not accurate in practice since data centres typically contain multiple-generations of physical devices with different processing capabilities~\cite{google_study_het}. Processing speeds of servers can vary also due to the presence of various types of acceleration devices such as GPUs, FPGAs and ASICs~\cite{GPU_het,FPGA_het}. Finding the delay optimal scheme for such systems remains an open problem to this date. 
For such systems, speed-unaware schemes such as JSQ and SQ($d$), designed primarily for homogeneous systems,
may perform very poorly~\cite{gardiner_perf,arpan_tcns,bramson2012asymptotic}.
Speed-aware schemes which assign jobs based on both queue lengths and server speeds can significantly outperform speed-unaware schemes. To see this, consider a speed-aware version of JSQ, henceforth referred to as the {\em Speed-Aware JSQ (SA-JSQ)} scheme, in which each arrival is assigned to a server with the highest speed  among the ones with the minimum queue length. In Figure~\ref{fig:1}, we compare the JSQ scheme with the SA-JSQ scheme for a system with two types of servers. Inter-arrival and service times of jobs are assumed to be exponentially distributed with rates $N\lambda$ and $1$, respectively. We observe that the average response time of jobs is almost $60$\% lower under the SA-JSQ scheme than under the JSQ scheme for loads ($\lambda$) less than $0.5$. This example clearly highlights the need to design speed-aware schemes which can lead to much improved performance compared to speed-unaware schemes for heterogeneous systems. 

There are various ways in which server speeds can be incorporated into job dispatching decisions. A natural scheme is to send each job to the queue with the shortest expected delay (SED), i.e., the queue where the ratio of the queue length to the server-speed is minimum. As we shall see later, for finite $N$ this scheme performs slightly better than the SA-JSQ scheme discussed above but this difference vanishes quickly as the number of servers increases. Thus, the key questions in the context of heterogeneous systems are the following: {\em what is the minimum achievable mean response time of jobs in heterogeneous systems?} and {\em is it possible achieve this minimum with some simple speed-aware variants of the JSQ scheme in the limit as the system size becomes large?}

\subsection{Contributions}

To answer the questions above we need to find a lower bound on the mean-response time of jobs in heterogeneous systems and analytically characterise the performance of speed-aware
dispatching schemes in the limit as $N \to \infty$. These are the primary aims of this paper. 
Our key finding is that in the fluid limit the SA-JSQ scheme achieves the lowest
possible mean response time of jobs for heterogeneous systems.
Note that to prove such optimality results, it is essential to develop the necessary tools to analyse job assignment schemes for heterogeneous systems. One would expect the fluid limit techniques applicable to homogeneous systems would easily generalise to heterogeneous systems. However, this turns out to be not the case. As explained later, the stochastic coupling technique, used to obtain uniform bounds on the stationary queue lengths in homogeneous systems, is not generalizable to heterogeneous systems. Consequently, we develop a new method based on exponential Lyapunov functions to obtain similar uniform bounds for heterogeneous systems. This method is sufficiently general to be applicable to other systems where a direct coupling is difficult to construct. Furthermore, the increased complexity of the underlying state space for heterogeneous systems results in a fluid limit which depends on the invariant distribution of a multi-dimensional birth-death process. Deriving exact analytical expressions for this invariant distribution is difficult. We find properties that are sufficient to show the existence, uniqueness and global stability of the fluid limit. From this point onward, we use the term `optimality' to refer to optimality in the fluid limit unless otherwise mentioned. Our contributions are listed below:

\begin{enumerate}
    \item {\em Stochastic comparison and lower bound}: Our first contribution is to compare the heterogeneous system, which consists of separate queues, with a similar heterogeneous system where all the servers serve a single central queue. We show that the system with separate queues stochastically dominates the system with the central queue. From this stochastic comparison result and the fluid limit of the system with central queue, we obtain a lower bound on the mean response time of jobs that holds for any dispatching scheme applied to the system with separate queues. It is important to note that studying fluid limit of the system with a single central queue and heterogeneous servers can be of independent interest. 
    
    \item {\em Stability and tightness}: We next show that the heterogeneous system is stable under the SA-JSQ policy and obtain uniform bounds (not depending on the system size) on the tails of the stationary queue length distribution. The latter result is required to show the tightness of the sequence of stationary distributions indexed by the system size. The usual approach for proving tightness via coupling does not work in the heterogeneous setting as it is difficult to construct a coupling that maintains the desired stochastic  dominance. Thus, we take a new approach and establish this result by analysing the drift of an exponential Lyapunov function. This technique is sufficiently general to be applicable to other systems where direct stochastic comparison is difficult. 
    
    \item {\em Fluid limit analysis}: Our final contribution is the fluid limit analysis of the SA-JSQ scheme in the heterogeneous system. This analysis turns out to be considerably more challenging than conventional fluid analysis due to the increased complexity of the underlying state space for heterogeneous systems and the separation of two time-scales in the fluid limit. The main idea in establishing the fluid limit of SA-JSQ is to first use the martingale functional central limit theorem (FCLT) and then characterise the limit of any convergent subsequence using Lemma~2 and Theorem~3 of~\cite{Hunt1994}.
    The fluid limit in the heterogeneous setting turns out to be a function of the invariant distribution of a multi-dimensional Markov process whose exact analysis seems intractable. However, we are able to prove properties of the invariant distribution sufficient to characterise the fluid limit and its unique fixed point. Using these properties and the monotonicity of the system we also show that the fluid limit is globally stable. This final result establishes the asymptotic optimality of the SA-JSQ policy for heterogeneous systems.
\end{enumerate}

\subsection{Related works}

There exists a vast literature on load-balancing policies for multi-server systems. 
Here, we only discuss the works that are most relevant to our paper. For a comprehensive review of existing works, we refer the reader to~\cite{van2018scalable}.

The JSQ policy was first shown to minimise the average delay of jobs for finite systems consisting of identical servers in~\cite{Winston1977optimality} under the assumption of Poisson arrivals and exponential service times. This optimality result was later extended to general stochastic arrival processes and service-time distributions
with non-decreasing hazard rates in~\cite{Weber1978}, to queues with state-dependent service rates in~\cite{johri1989optimality}, and to systems with finite buffer-buffers and general batch arrivals in~\cite{hordijk1990optimality}. 
Recent works~\cite{Mukherjee2018,Eschenfeldt2018} have considered the fluid and diffusion limits of the JSQ scheme. In the fluid limit, it has been shown in~\cite{Mukherjee2018} that the 
fraction of servers with two or more jobs converges to zero under the JSQ scheme. 
This implies that in the fluid limit all jobs 
find an idle server to join.
In the Halfin-Whitt regime, where the normalized arrival rate $\lambda$ varies with the system size $N$ as $\lambda=1-\beta/\sqrt{N}$ for some $\beta >0$,
it has been shown in~\cite{Eschenfeldt2018}
that the diffusion-scaled process
approaches to a two-dimensional reflected Ornstein-Ulhenbeck (OU) process as $N \to \infty$. The stationary distribution of this OU process has been studied in~\cite{Braverman2020} which establishes that the steady-state fraction of servers with exactly two jobs scales as $O(1/\sqrt{N})$ and the fraction of servers with more than three jobs scales as $O(1/N)$.


Another line of works explores load balancing schemes which require less communication between the servers and the job dispatcher. The SQ($d$) scheme in which each arrival is assigned to the  shortest among $d$ randomly sampled queues was first analyzed independently in~\cite{Vvedenskaya1996} and \cite{Mitzenmacherthesis}.
Using mean-field analysis, it was shown that for $d \geq 2$ the stationary queue length distribution has a super-exponentially decaying tail for large system sizes. Thus, by querying only $d \geq 2$ servers at every arrival instant a significant reduction in the average delay can be obtained in comparison to $d=1$.
The Join-the-Idle-Queue (JIQ) further reduces
communication overhead by keeping track of only the idle servers in the system. At each arrival instant it sends the incoming job to an idle server if one is available; otherwise, the job is sent to a randomly sampled server.
This scheme was first proposed and analyzed by Lu {\em et al.} in~\cite{Lu_JIQ_2011} and it was shown that in the fluid limit the JIQ scheme achieves the same performance as the JSQ scheme. Thus, the JIQ scheme is asymptotically optimal for homogeneous systems in the fluid limit.

Relatively few works consider load balancing
in heterogeneous systems. The SQ($d$) scheme for
heterogeneous systems has been analyzed in
\cite{shroff_heavy_traffic_2017,arpan_tcns,Makowski_SQd_2014}. While~\cite{Makowski_SQd_2014} and
\cite{shroff_heavy_traffic_2017} analyze
the performance of the SQ($d$) policy in
light and heavy traffic regimes for finite system sizes, respectively,~\cite{arpan_tcns} considers its performance in the mean-field regime. It has been shown that the SQ($d$) scheme suffers from a reduced stability region in heterogeneous systems due to infrequent sampling of faster servers. Subsequent works~\cite{gardiner_perf,arpan_ssy} have studied variations of the SQ($d$) scheme, aimed at improving its performance in heterogeneous systems while retaining the maximal stability region.
Heterogeneity-aware load balancing (HALO) schemes such Random, Round-Robin (RR), and SQ(d) have been analyzed in~\cite{Gandhi2015halo} for heterogeneous processor sharing systems. In particular, the optimal load splitting  has been obtained for the Random scheme and then used for the other schemes. It has been shown that these schemes result in good performance for all system sizes.
The JIQ scheme has been analyzed in the heterogeneous setting by Stolyar~\cite{Stolyar2015}; it has been shown that the average waiting time of jobs under the JIQ scheme approaches to zero in the fluid limit. 
From this result, it follows intuitively that a scheme which further distinguishes between idle servers by their speeds can only result in a smaller processing time of jobs. However, to establish this formally, it is necessary to construct a coupling that maintains the dominance 
of speed-aware schemes over speed-unaware schemes for sufficiently large $N$.
Constructing such coupling in turn requires establishing the fluid limit for speed-aware schemes, which brings us back to the problem considered in this paper.


Some recent works, e.g., \cite{debankur_constrained_2021,Weng2020}, study load balancing schemes for systems where jobs are constrained to be served only by specific subsets of servers. In these works, the focus is on finding conditions on the compatibility constraints  such that 
the performance of classical load balancing algorithms such as JSQ, JIQ and SQ($d$) remain asymptotically the same as in systems without compatibility constraints. 
In the work by Weng {\em et al.}~\cite{Weng2020}, a scheme
similar to the SA-JSQ scheme has been considered for constrained heterogeneous systems with finite buffer sizes. They prove the asymptotic optimality of this scheme using Lyapunov drifts methods.
Although the scheme is technically the same as the SA-JSQ scheme for constrained systems, their analysis crucially relies on the assumption that the queues have finite buffer sizes. Hence, in their setting tightness and stability results follow immediately. In contrast, in the setting considered in this paper, the queues have infinite buffer sizes and proving optimality in this setting requires showing the tightness of the stationary distributions. As discussed before, this is one of the key difficulties in our model and we employ a new technique to establish tightness. Furthermore, the drift technique, applicable to finite-buffer systems, is difficult to generalise to our setting. Instead, to prove asymptotic optimality, we use martingale methods outlined in~\cite{Whitt2007}.

In this paper, the drift of the fluid-scaled Markov chain depends on the un-scaled  state of the system through an indicator function.
This dependence ultimately gives rise to a time-scale separation in the fluid limit. A similar drift structure involving indicator functions was studied in~\cite{gast2012markov}. Here, the challenge is to tackle the discontinuity of the drift at the boundaries of the state-space. It has been shown in~\cite{gast2012markov} that in such cases the finite system can be viewed as a stochastic approximation of a differential inclusion. 
In~\cite{ayesta2012scheduling}, the author study optimal scheduling schemes for users competing for a common resource.
They show that a tie-breaking rule similar to JSQ, where users with highest departure probability are preferred, leads to asymptotic optimality. 

\subsection{General Notations}
We use the following notations throughout the paper. The sets $\mb{R,N, Z_+}$ denote the set of 
real numbers, natural numbers, and non-negative integers, respectively.
We also use the set $\bar \Z_+$ to denote the extended set $\Bar{\Z}_+=\Z_+ \bigcup \{\infty\}$.
For $x,y \in \mb{R}$, we use $x\wedge y$, $x \vee y$, and $(x)_+$ to 
denote $\max(x,y)$, $\min(x,y)$, and $\max(x,0)$, respectively. 
For any $n \in \mb{N}$, $\sbrac{n}$ denotes the set $\cbrac{1,2,\ldots,n}$. For any complete separable metric space $E$, we denote $D_{E}[0,\infty)$ to be the set of all {\em cadlag} functions from $[0,\infty)$ to $E$ endowed with the Skorohod topology. Moreover, the notation $\mathscr{B}(E)$ is used to denote the Borel sigma algebra generated by the set $E$. 
The notation $\Rightarrow$ is used for weak convergence. We use $\indic{A}$ to denote
indicator function for set $A$.

\subsection{Organisation}
The rest of the paper is organized as follows.
In Section~\ref{Sec:system_model}, we introduce the system model and define the SA-JSQ policy. In Section~\ref{Sec:RP}, we obtain a lower bound on the mean response time of jobs that holds for any scheme in the heterogeneous system by comparing the system with a similar system having a central queue. In Section~\ref{sec:mf_sajsq}, we state
our main results for the SA-JSQ policy.
We prove the stability of the SA-JSQ scheme in Section~\ref{sec:Stability_SA-JSQ} and the uniform bounds on the tails of the stationary queue lengths distribution in Section~\ref{sec:uniform_bounds}. The monotonicity property for SA-JSQ for a finite $N$ system is shown in Section~\ref{sec:monotone}. In Section~\ref{subsec:Process_Level_Convergence_SA-JSQ}, we prove the process convergence result for SA-JSQ and characterise its fixed point in Section~\ref{sec:Global_stability_Limit_interchange}. The proof of resource pooled optimality is given in Section~\ref{sec:Rp_optimality}. Numerical studies comparing different schemes under heterogeneous setting is given in Section~\ref{sec:Numerical_Studies}. Finally, we conclude the paper in Section~\ref{sec:Conclusion}.

\section{System Model} 
\label{Sec:system_model}
We consider a system $\mc M_N$ consisting of $N$ parallel servers, 
each with its own queue of infinite buffer size.
The servers are assumed to be {\em heterogeneous}, i.e., they have different service rates. Specifically, we 
assume that there are $M$ different server {\em types} or {\em pools}. 
Each type $j \in [M]$ server has a service rate of $\mu_j$. The proportion of type $j \in [M]$ servers in the system is assumed to be fixed at $\gamma_j \in [0,1]$ with $\sum_{j\in [M]}\gamma_j=1$.
We further assume without loss of generality that $\mu_1 > \mu_2 > \ldots > \mu_M$
and $\sum_{j \in [M]} \gamma_j \mu_j=1$ (normalised system capacity is unity). For simplicity of exposition we also assume that there exists $N \in \mb N$ such that $N\gamma_j \in \mb{N}$ for all $j \in [M]$.\footnote{Our asymptotic results do not depend on this assumption. However,
the results for finite systems need to be modified slightly if this is not the case.}


Jobs are assumed to arrive at the system according to a Poisson process with a rate $N\lambda$. 
Each job requires a random amount of work, independent and 
exponentially distributed with mean $1$. The inter-arrival and service 
times are assumed to be independent of each other. 
Upon arrival, each job is assigned to a server where it either immediately
receives service (if the server is idle at that instant) or joins the
corresponding queue to be served later. 
The queues are served according to the First-Come-First-Server (FCFS) scheduling discipline. 
For this system, the term {\em queue length} will refer to the 
total number of jobs in a queue including the current job in service.


Our main interest is to find a job assignment policy that 
minimises the steady state mean response time of jobs in the system. 
To this end, we shall analyse
a modified version of the classical Join-the-Shortest-Queue (JSQ) policy which is known to be optimal
for homogeneous systems~\cite{Weber1978}.
We shall refer to the modified policy as the {\em Speed-Aware JSQ} or the {\em SA-JSQ} policy. It is defined as follows

\begin{defn}
Under the SA-JSQ policy, upon arrival of a job, 
it is sent to a server with the minimum queue length 
among all the servers in the system.
Ties among servers of different types are broken by choosing the server type with the maximum speed and ties among servers
of the same type are broken uniformly at random.
\end{defn}

\begin{figure}[h!]
  \centering
  \includegraphics[width=8cm]{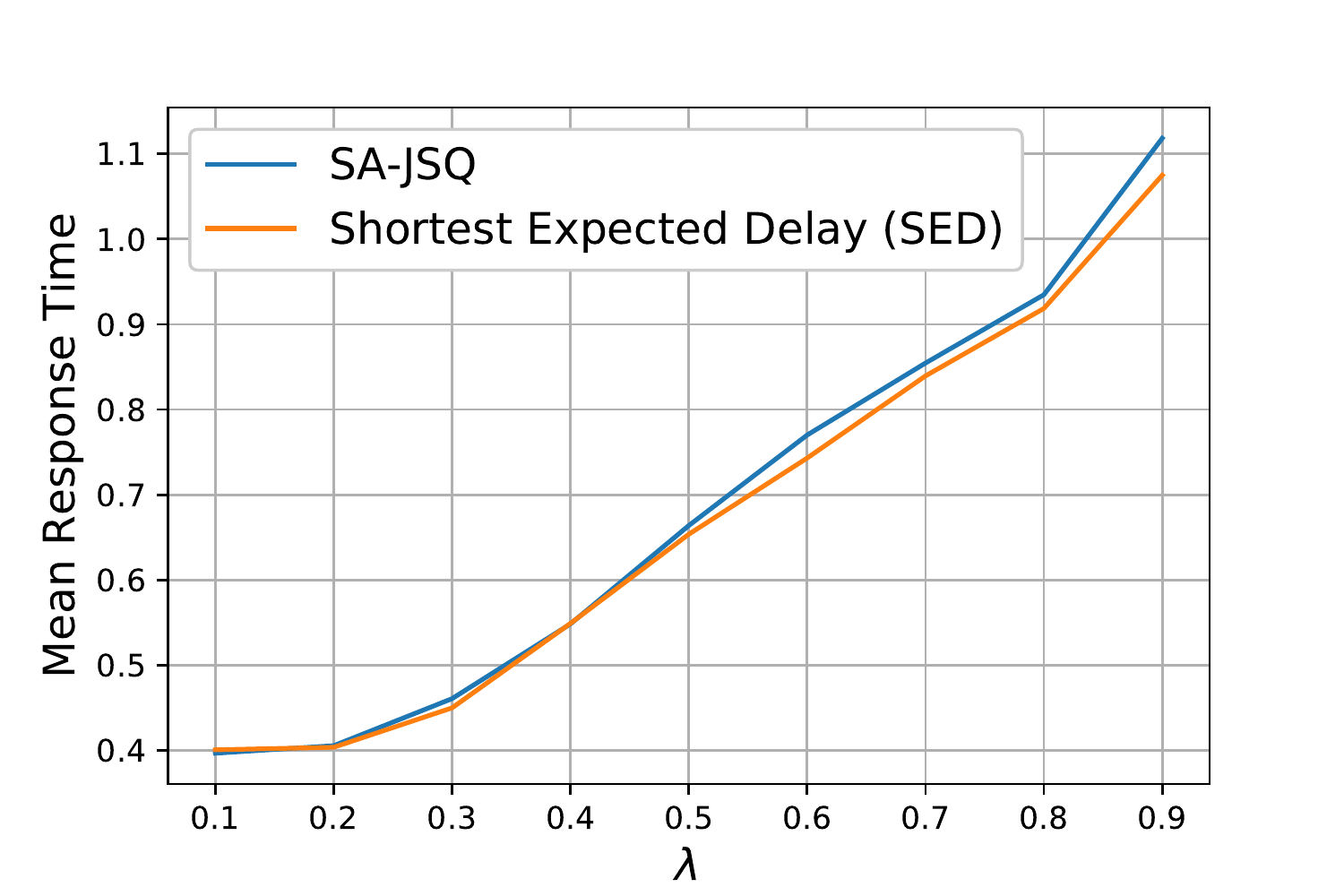}
\caption{Comparison of mean response time of jobs between SED and SA-JSQ schemes under heterogeneous system with $N=50$.  We choose
$\gamma_1=1-\gamma_2=1/5$, and $\mu_1=4\mu_2=20/8$.}
\label{fig:finite_n}
\end{figure}

Note that unlike the classical JSQ policy, which breaks ties uniformly at random, 
the SA-JSQ policy breaks ties among servers of different types by choosing the server type with the
maximum service rate.
Furthermore, unlike the classical JSQ policy,
the SA-JSQ policy is not optimal
for the heterogeneous system for finite
$N$: numerical simulations shown in Figure~\ref{fig:finite_n} show that
a scheme which assigns jobs
to servers based on the shortest
expected delay (SED) (i.e., queue length divided by
the service rate) performs marginally better than
the SA-JSQ policy for small system sizes.
However, as we show later in the paper, the gap vanishes as $N \to \infty$. More specifically, we show that
in the limit as $N \to \infty$ there is no better scheme than SA-JSQ for heterogeneous systems.

\section{Lower Bound on the Mean Response Time}
\label{Sec:RP}

A key ingredient in establishing the optimality of the SA-JSQ policy
is finding a uniform lower bound 
on the steady state mean response time 
of jobs in system $\mc M_N$
for all $N$ under any stationary 
job assignment policy $\Pi$.
We later show
that this lower bound is achieved by the SA-JSQ policy
as $N \to \infty$.

To obtain this lower bound, we show that the 
total number of jobs in system $\mc M_N$ under any stationary policy $\Pi$ 
stochastically dominates the total number of jobs
in a similar heterogeneous system $\mc M_N'$
working under a specific policy, referred to as
the {\em Join-the-Fastest-Free-Server (JFFS)} policy.


{\em Description of system $\mc M_N'$}:
In system $\mc M_N'$, the servers remain the same as in system $\mc M_N$, i.e., there are $N$ servers in total
with $N\gamma_j$ of them having service rate $\mu_j$
for $j \in [M]$.
However, unlike system $\mc M_N$, 
we have a central queue for the whole system. 
Upon arrival, a job joins the central queue if all the servers are busy. 
Otherwise, the job is assigned to an idle server.
Here, the central queue consists of only those jobs that are waiting for service
and it excludes the jobs already in service.

{\em JFFS policy}:
If multiple servers are idle and the queue is non-empty,
a decision needs to be made on which server the next job should be assigned to. 
To this end, we consider the following job assignment 
policy referred to as the {\em Join-the-Fastest-Free-Server (JFFS)} policy.

\begin{defn}
Under the JFFS policy, if idle servers are available and there are jobs waiting in the central 
queue, then the head-of-the-line (HOL) job is assigned to the idle server with the highest service rate. 
Ties within the same server type are broken uniformly at random. 
\end{defn}

Hence, unlike system $\mc M_N$ where each server only serves its own queue, 
in system $\mc M_N'$, the central queue is served by all the servers. 
This system is, therefore, referred to as the {\em resource pooled system}~\cite{Song2015}. 


The evolution of the resource pooled system $\mc M_N'$
under the JFFS policy
can be described by the Markov chain 
$Z^{(N)}=(Z^{(N)}(t),t\geq 0)$, where $Z^{(N)}(t) \in \mb{Z}_+$
denotes the total number of jobs in system $\mc M_N'$ 
at time $t$ under the JFFS policy. 
Clearly, under the JFFS policy,
the number of busy type $j$
servers at time $t \geq 0$
is given by $(Z^{(N)}(t)-\sum_{i=1}^{j-1}N\gamma_i)_+\land N\gamma_j$.
Hence, the transition rate $q^{Z^{(N)}}(k,l)$
of $Z^{(N)}$ from state $k$ to state 
$l$ is given by 
\begin{equation}
    q^{Z^{(N)}}(k,l)=\begin{cases}
                            N\lambda, &\text{ if } l=k+1,\\
                            \sum_{j \in [M]}\mu_j((k-\sum_{i=1}^{j-1}N\gamma_i)_+\land N\gamma_j), &\text{ if } l=k-1,\\
                            0, &\text{ otherwise}.
    \end{cases}
\end{equation}
The following two propositions
characterise the stationary behaviour of the chain $Z^{(N)}$
for finite $N$ and for $N \to \infty$, respectively. 


\begin{proposition}
\label{thm:stability_RP}
For $\lambda < 1$, the Markov chain  $Z^{(N)}$ is 
positive recurrent. Furthermore,
if $Z^{(N)}(\infty)$ denotes the stationary number of jobs in
the system, then 

\begin{equation}
\label{eqn:sup_bound}
    \sup_{N} \frac{\E\sbrac{Z^{(N)}(\infty)}}{N}\leq \lambda +\frac{\lambda}{1-\lambda}.
\end{equation}
\end{proposition}

\begin{proposition}
\label{thm:mf_rp}
If $\lambda <1$, then 
for system $\mc M_N'$ operating under the JFFS
policy we have as $N \to \infty$

\begin{equation}
 \frac{{Z^{(N)}(\infty)}}{N} \Rightarrow z^*\define\max_{j\in [M]}\brac{
 \sum_{i=1}^{j-1}\gamma_i+ \frac{\lambda - \sum_{i=1}^{j-1}\mu_i\gamma_i}{\mu_j}}
 \label{eq:def_q*}
\end{equation}
\end{proposition}

Proposition~\ref{thm:stability_RP} implies that the 
resource pooled system is
stable for $\lambda < 1$.
It also provides a uniform
upper bound on the steady-state 
expected (scaled) number of jobs in
the system for all $N$.
These results are utilised to prove Proposition~\ref{thm:mf_rp}
which shows that the steady-state (scaled)
number of jobs in the resource pooled system converges weakly to $z^*$ as $N \to \infty$, where $z^*$ is
as defined in~\eqref{eq:def_q*}.

In~\eqref{eq:def_q*}, 
the expression within the brackets 
on the RHS
is maximised for 
$j=j^* \in [M]$ iff $\lambda \in \left[\sum_{i=1}^{j^*-1} \gamma_i \mu_i,\sum_{i=1}^{j^*} \gamma_i \mu_i\right)$.
Note that the existence of such a $j^*$ is guaranteed
because $\lambda < 1=\sum_{j \in [M]} \gamma_j \mu_j$.
For $j^*$ as defined above, it is easy to see that the fraction of busy type $j$ servers is one for all $j \in [j^*-1]$ and zero for all $j \geq j^*+1$. Furthermore, in pool $j^*$, the fraction of busy servers is given by $\frac{\brac{\lambda-\sum_{i=1}^{j^*-1}\gamma_i \mu_i}_+}{\gamma_{j^*} \mu_{j^*}}.$ Since $z^* < 1$, the central queue in the limiting
system is empty in the steady-state.


To show that the system $\mc M_N$
operating under any stationary policy $\Pi$ 
stochastically dominates the system $\mc M_N'$ 
operating under the JFFS policy,
we describe the state of the system $\mc M_N$ under a stationary policy $\Pi$ by the 
Markov chain $\mf X^{(N,\Pi)}=(X^{(N,\Pi)}_{i,j}(t), t \geq 0)$, where $X^{(N,\Pi)}_{i,j}(t)$ denotes the number of type $j$ servers with at least $i$ jobs. Let $R^{(N,\Pi)}(t)=\sum_{i,j} X^{(N,\Pi)}_{i,j}(t)$ denote the total number of jobs in system $\mc M_N$ 
at time $t \geq 0$ under policy $\Pi$. 
We have the following stochastic comparison result.

\begin{thm}
\label{thm:RP_optimality}
For any stationary policy $\Pi$ if $Z^{(N)}(0)\leq R^{(N,\Pi)}(0)$, then 
the processes $Z^{(N)}$ and $\mf X^{(N,\Pi)}$ can be constructed on the same
probability space 
such that $Z^{(N)}(t)\leq R^{(N,\Pi)}(t)$ for all $t\geq0$, i.e., $Z^{(N)}(t)\leq_{st} R^{(N,\Pi)}(t)$, where $\leq_{st}$ implies stochastic dominance.
\end{thm}

The above theorem implies that
if both $Z^{(N)}$ and $R^{(N,\Pi)}$
are ergodic chains,  then $Z^{(N)}(\infty) \leq_{st} R^{(N,\Pi)}(\infty)$, where $Z^{(N)}(\infty)=\lim_{t \to \infty} Z^{(N)}(t)$ and $R^{(N,\Pi)}(\infty)=\lim_{t\to \infty}R^{(N,\Pi)}(t)$ denote the stationary number of jobs in $\mc M_N'$ and $\mc M_N$, respectively.
This further implies that 
$\expect{Z^{(N)}(\infty)} \leq \expect{R^{(N,\Pi)}(\infty)}$ which by Little's law 
gives the desired lower bound on the stationary 
mean response time of jobs in $\mc M_N$ under any policy $\Pi$.

By Little's law, the steady-state mean response time, $\bar T_N'$,
of jobs in $\mc M_N'$ under the JFFS policy
is given by $\bar T_N'=\expect{Z^{(N)}(\infty)}/N\lambda$.
Hence, Proposition~\ref{thm:mf_rp} implies 
$\lim_{N \to \infty} \bar T_N'=\frac{z^*}{\lambda}$.
Combining the above result with the stochastic comparison
result in Theorem~\ref{thm:RP_optimality}, we can conclude the following.

\begin{cor}
\label{cor:optimality}
The steady-state mean response time, 
$\bar T_{N,\Pi}$, of jobs
in system $\mc M_N$ under any stationary policy $\Pi$
satisfies
\begin{equation}
    \label{eq:lower_bound}
    \liminf_{N \to \infty} \bar T_{N,\Pi}\geq \frac{z^*}{\lambda}
\end{equation}
\end{cor}

In subsequent sections, we establish that
the above lower bound is achieved with equality
when SA-JSQ is employed as the job assignment policy in
$\mc M_N$. This will imply the asymptotic optimality of
the SA-JSQ policy for $\mc M_N$.
\section{Analysis setup and main results}
\label{sec:mf_sajsq}

In this section, we introduce the key
ingredients and notations for our analysis of the SA-JSQ policy. In particular, we discuss the different state-descriptors used in the analysis and discuss
some key properties of the space they belong to. We also state our main results for the SA-JSQ policy and discuss their implications.

The state
of the system $\mc M_N$ at any time $t \geq 0$ 
under the SA-JSQ policy can be described in two different ways. These are as defined below:

\begin{enumerate}
    \item {\em Queue-length descriptor}: We define the queue-length vector
    at time $t\geq 0$ as
    $$\mf Q^{(N)}(t)=(Q_{k,j}^{(N)}(t),k \in [N\gamma_j], j \in [M]),$$ 
    where $Q_{k,j}^{(N)}(t)$ denotes
    the queue length of the $k^{\textrm{th}}$ server
    of type $j$ at time $t$.
    
    \item {\em Empirical measure descriptor}: We define the empirical
    tail measure on the queue lengths at time $t$ as 
    $$\mathbf{x}^{(N)}(t)=({x}^{(N)}_{i,j}(t),i \geq 1, j\in[M]),$$
    where ${x}^{(N)}_{i,j}(t)$ denotes the fraction of type $j$ servers with at least 
    $i$ jobs at time $t$. For completeness, we set
    $x^{(N)}_{0,j}(t)=1$ for all $j \in [M]$ and all $t\geq 0$.
\end{enumerate}
It follows from the Poisson arrival and exponential
job size assumptions
that both processes 
$\mf Q^{(N)}=(\mf Q^{(N)}(t), t\geq 0)$ and $\mathbf{x}^{(N)}=(\mathbf{x}^{(N)}(t), t\geq0)$ are Markov. It is possible to switch from first descriptor
to the second by noting the following 
\eqn{
\label{eqn:Realtion_Q_X}
x_{i,j}^{(N)}(t)=\frac{1}{N\gamma_j}\sum_{k\in[{N\gamma_j}]}\indic{Q_{k,j}^{(N)}(t)\geq i}, \ \forall i\geq1, \ j\in[M].
}
We use 
both descriptors above to state and prove our results.
%

One important observation to make is that
the SA-JSQ policy does not distinguish between
servers of the same type. Hence, if
the distribution $\nu$ of the initial state of the system is {\em exchangeable}
within each server type, i.e., 
if for each $j$ we have 
$\mb{P}_\nu(Q_{k_l,j}(0)=r_l,l \in [r])=
\mb{P}_\nu(Q_{\sigma(k_l),j}(0)=r_l, l \in [r])$
for all permutations $\sigma$
of the numbers $k_1,k_2,\ldots,k_r \in [N\gamma_j]$,
then for any time $t \geq 0$, we must have
$\mb{P}_\nu(Q_{k_l,j}(t)=r_l,l \in [r])=
\mb{P}_\nu(Q_{\sigma(k_l),j}(t)=r_l, l \in [r])$, where $\mb{P}_\nu$ is the probability
measure conditional on the initial state being distributed according to $\nu$.
This property is crucial in connecting the distribution
of $\mf Q^{(N)}$ to that
of $\mf x^{(N)}$. In particular, by taking
expectation on both sides of \eqref{eqn:Realtion_Q_X}
we have
\begin{align}
    \E_\nu[x_{i,j}^{(N)}(t)]=\frac{1}{N\gamma_j}\sum_{k\in [N \gamma_j]} \mb{P}_\nu(Q_{k,j}(t) \geq i)
    =\mb{P}_\nu(Q_{k,j}(t) \geq i),
\end{align}
where in the second equality we employ the exhangeability property discussed above.
Clearly, the process $\mf Q^{(N)}$ takes values in $\mb{Z}_+^N$ and
the process $\mf x^{(N)}$ takes values in the space ${S}^{(N)}$
defined as 
\begin{equation}
\label{eqn:finite_N_Space}
    {S}^{(N)}\define \{ \mf s=(s_{i,j}):N\gamma_js_{i,j}\in \Z_+, 1\geq s_{i,j} \geq s_{i+1,j} \geq 0 \ \forall  i\geq 1,j\in[M]\}.
\end{equation}
Note that for finite $N$, the space ${S}^{(N)}$ is countable
since each $x_{i,j}^{(N)}$ can only take finitely many values. 
We further define the space $S$ 
as follows
\begin{equation}
\label{eqn:finite_infinite_Space}
\begin{aligned}
     S\define& \{ \mf s=(s_{i,j}):1\geq s_{i,j} \geq s_{i+1,j} \geq 0, \forall  i\geq 1,j\in[M], \norm{\mf s}_1<\infty\},
\end{aligned}
\end{equation}
where the $\ell_1$-norm, denoted by $\norm{\cdot}_1$, is defined as $\norm{\mf s}_1\define\max_{j \in [M]} \sum_{i\geq 1} |s_{i,j}|$
%
%
for any $\mf s \in S$. 
It is easy to verify that the space $S$ 
is complete and separable under the $\ell_1$-norm.
Furthermore, if the system's state $\mf x^{(N)}$ belongs to
$S^{(N)} \cap S$, then there are finitely
many jobs in the system. 
Starting with a state in $S^{(N)} \cap S$ we can ensure that chain $\mf x^{(N)}$ remains in $S\cap S^{(N)}$ for all $t\geq 0$
only if the process $\mf x^{(N)}$, or equivalently the process $\mf Q^{(N)}$, 
is positive recurrent.
Our first main result states that this is indeed the case when
$\lambda < 1$.

\begin{thm}
\label{thm:SA-JSQ_Stability}
The process $\mf Q^{(N)}$ or equivalently $\mf x^{(N)}$ is positive recurrent for each $\lambda < 1$ and each $N$.
\end{thm}

The theorem above implies that for $\lambda < 1$ the stationary distributions
of $\mf Q^{(N)}$ and $\mf x^{(N)}$ exist
and they are unique.
Let $\mathbf{Q}^{(N)}(\infty)=\lim_{t \to \infty}\mathbf{Q}^{(N)}(t)$ (resp.
$\mf x^{(N)}(\infty) =\lim_{t \to \infty} \mf x^{(N)}(t)$) denote 
the state of the system distributed according to the stationary distribution of $\mf Q^{(N)}$ (resp. $\mf x^{(N)}$). 
The stability of the system
guarantees that 
the steady-state expected number of
jobs in the system is finite, i.e.,
$\mf x^{(N)}(\infty)$ belongs to 
$S$ almost surely.
However, to establish the asymptotic
optimality of the SA-JSQ policy, we 
require the stronger result that
the sequence $(\mf x^{(N)}(\infty))_N$
of stationary states indexed by the system size is tight in $S$. Here, it is important to observe that $S$ is not a compact space. The relatively compact subsets of $S$ and the tightness criterion
for sequences in $S$ are stated in~\ref{sec:compact_sets}. The tightness criterion essentially requires 
the expected tail sums
$\E[\sum_{j \in [M]} \sum_{i \geq l} x_{i,j}^{(N)}(\infty)]$ approach to zero
as $l \to \infty$ uniformly in $N$.
To show that $(\mf x^{(N)}(\infty))_N$ satisfies this criterion, we need the following important result
which states that the stationary tail probabilities
of the queue lengths decay exponentially
and the decay rate is uniform in $N$. 

\begin{thm}
\label{thm:tail_bound}
If $\lambda < 1$, then under the SA-JSQ scheme, for each $j \in [M]$, each $k \in [N\gamma_j]$, and each $l \geq 1$ the following bound holds for all $\theta \in [0, -\log \lambda)$ 
\eqn{
\sup_N \mathbb{P}(Q_{k,j}^{(N)}(\infty)\geq l)\leq C_j(\lambda, \theta)e^{-l \theta},
\label{eq:tail_bound}
}
where $C_j(\lambda,\theta)=(1-\lambda)/(\mu_j \gamma_j (1-\lambda e^\theta)) > 0$.
\end{thm}

Our next main result establishes a key monotonicity property of the system for finite $N$. For two states $\mf q$
and $\mf{\tilde{q}}$ in $\mb Z_+^N$
we say $\mf q \leq \mf{\tilde{q}}$
if $q_{k,j} \leq \tilde q_{k,j}$
for each $k\in [N \gamma_j],j \in [M]$.
The inequality $\mf s \leq \mf{\tilde s}$ is similarly defined for $\mf s,\mf{\tilde s} \in S$. We have the following result.

\begin{thm}
\label{thm:SA-JSQ_monotone}
Consider two systems with initial states $\mf Q^{(N)}(0)$ and $\tilde{\mf Q}^{(N)}(0)$ (resp. $\mf x^N(0)$ and $\tilde{\mf x}^N(0)$) satisfying
$\mathbf{Q}^{(N)}(0)\leq\tilde{\mathbf{Q}}^{(N)}(0)$ (resp. $\mf x^N(0) \leq \tilde{\mf x}^N(0)$).
Let $\mf Q^{(N)}$ and $\tilde{\mf Q}^{(N)}$
(resp. $\mf x^N$ and $\tilde{\mf x}^N$)
denote the corresponding processes describing
the two systems under the SA-JSQ policy.
Then, the
processes $\mf Q^{(N)}$ and $\tilde{\mf Q}^{(N)}$ (resp. $\mf x^N$ and $\tilde{\mf x}^N$) can be constructed on the same probability space
such that $\mathbf{Q}^{(N)}(t)\leq \tilde{\mathbf{Q}}^{(N)}(t)$ (resp. $\mf x^N(t)\leq \tilde{\mf x}^N(t)$)
for all $t \geq 0$.
\end{thm}

The monotonicity
property stated above 
implies that if two systems,
both working under the SA-JSQ policy, start 
at two different initial states 
such that the queue lengths in the first system
are dominated by the corresponding queue lengths in the second system, then this dominance is maintained for all $t \geq 0$.
This property turns out to be
essential in establishing the global asymptotic stability of the fluid limit process.

Our next set of results characterise the asymptotic properties of the process $\mf x^{(N)}$ as $N \to \infty$. The first result states that the sequence $(\mf x^{(N)})_N$ of processes indexed by $N$ converges weakly
to a deterministic process $\mf x=(\mf x(t), t\geq 0)$ defined on $S$. To describe the limiting process,
we define $l_j(\mf s)=\min\{i:s_{i+1,j}<1\}$ for any state $\mf s\in S$  to be the minimum queue length in pool $j \in [M]$ in state $\mf s$. 

\begin{thm}
\label{thm:SA-JSQ_Process_level_convergence}
(Process Convergence):  
Assume ${\mathbf{x}}^{(N)}(0)\in S\cap S^{(N)}$ for each $N$ and $ {\mathbf{x}}^{(N)}(0) \Rightarrow {\mathbf{x}}(0)\in S$ as $N \rightarrow \infty$. Then, the sequence $(\mathbf{x}^{(N)})_{N \geq 1}$ is relatively compact in $D_S[0,\infty)$
and any limit $\mf x=(\mf x(t)=(x_{i,j}(t), i \geq 1, j \in [M]), t \geq 0)$
of a convergent  sub-sequence 
of $(\mathbf{x}^{(N)})_{N \geq 1}$
satisfies the following set of equations for all $t\geq0$, $i\geq1$ and $j\in[M]$
\EQN{
\label{eqn:SA-JSQ_fluid_process}
{x}_{i,j}(t)=x_{i,j}(0)+ \frac{\lambda}{\gamma_j} \int_0^t p_{i-1,j}(\mathbf{x}(u))du -\int_0^t \mu_j (x_{i,j}(u)-x_{i+1,j}(u))du,
}
where 
$p_{i-1,j}(\mathbf{s}) \in [0,1]$ is uniquely determined 
for each state $\mf s \in S$. Furthermore,
$\mf p(\mathbf{s})=(p_{i-1,j}(\mathbf{s}),i\geq1,j\in[M])$ for each $\mf s \in S$
satisfies the following properties
\begin{condenum}
    \item $\sum_{j\in[M]}\sum_{i\geq1}p_{i-1,j}(\mathbf{s})=1$ for all $\mf s \in S$, \label{p2}
    
    \item $p_{i-1,j}(\mathbf{s})=0$ for all $i\geq l_j(\mf s)+2$ and for all $j\in[M]$,\label{p3}
    
    \item If $l_j(\mf s)>0$ for some $j\in[M]$, then $p_{i-1,j}(\mathbf{s})=0$ for all $1\leq i\leq l_j(\mf s)-1$.\label{p4}
    
    \item If $l_1(\mf s)=0$, then $p_{0,1}(\mathbf{s})=1$.\label{p5}
    
    \item For some $j\in\cbrac{2,\dots,M}$. If $l_{k}(\mf s)\geq1$ for all $k\in[j-1]$ and $l_{j}(\mf s)=0$, then $p_{0,j}(\mf s)=1$.\label{p6}
\end{condenum}
\end{thm}

In the theorem above, $p_{i-1,j}(\mf s)$
can be interpreted as the limiting probability of an incoming arrival being assigned to a type $j$ server with queue length $i-1$ when the system is in state $\mathbf{s} \in S$. With this interpretation, the properties of $p_{i,j}(\mf s)$ listed above follow intuitively from the assignment rule under the SA-JSQ policy. Indeed, properties P\ref{p5} and P\ref{p6} state that if 
in state $\mf s \in S$
pool $j \in [M]$ is the pool with 
the highest speed containing
idle servers, then
with probability $1$, incoming jobs are assigned to idle servers in pool $j$ in state $\mf s$. 
Similarly, properties P\ref{p3} and P\ref{p4} imply that
for any state $\mf s \in S$
and any pool $j \in [M]$ with minimum queue-length $l_j(\mf s)$, 
jobs can only be assigned to
queues with lengths $l_j(\mf s)-1$ and
$l_j(\mf s)$ in the limiting system.

The exact expressions for $p_{i,j}(\mf s)$ are complicated as they depend on the stationary probabilities of a multi-dimensional Markov chain (we explain this more in Section~\ref{subsec:Process_Level_Convergence_SA-JSQ}).
In Theorem~\ref{thm:SA-JSQ_Process_level_convergence}, we only list the properties essential to characterise the fixed point $\mf x^*=(x_{i,j}^*)$ of the fluid limit $\mf x$. In our final result stated below, we show that
the fixed point $\mf x^*$ is unique  and globally attractive.

\begin{thm}
\label{thm:SA-JSQ_global_stability_limit_interchange}
\begin{temp}
    \item For $\lambda < 1$, there exists a unique fixed point $\mathbf{x}^*=(x^*_{i,j},i\geq1, j\in[M]) \in S$ of the fluid limit $\mf x$ described by~\eqref{eqn:SA-JSQ_fluid_process}, i.e., if $\mf x(0)=\mf x^*$ then $\mf x(t)=\mf x^*$ for all $t\geq 0$. Furthermore, the fixed point $\mf x^*$ is given by
\EQN{
\label{eqn:SA-JSQ_Fixed_Point}
\begin{aligned}
 x^*_{1,j}=\Big(1 \land \frac{(\lambda-\sum_{k=1}^{j-1}\mu_k \gamma_k)_+}{\mu_j \gamma_j}\Big) \ \forall j\in [M], \ \ 
x^*_{i,j}=0 \ \ \forall \ i\geq2, j\in[M].
\end{aligned}
}\label{fixed_point}
    \item  (Global Stability): If $\lambda < 1$, then for any solution $\mf x$ of~\eqref{eqn:SA-JSQ_fluid_process} with
    $\mf x(0) \in S$ converges to $\mf x^*$ component-wise, i.e., $x_{i,j}(t) \to x_{i,j}^*$ as $t \to \infty$ for all $i\geq 1$ and for all $j\in [M]$.\label{global_stab_sa_jsq}
\item (Interchange of Limits): Let $\lambda < 1$. Then, the sequence $(\mf x^{(N)}(\infty))_N$ converges weakly to $\mf x^*$, i.e.,
$\mf x^{(N)}(\infty) \Rightarrow {\mathbf{x}^*}$ as $N \rightarrow \infty$.\label{limit_interchange}
\end{temp}
\end{thm}

Note that the last statement of the theorem implies that the sequence
of stationary distributions indexed by
the system size concentrates on the 
point $\mf x^*$. Furthermore, the first statement of the theorem implies
that in state $\mf x^* \in S$ the fraction of servers with two or more jobs is zero. This is similar to 
the classical JSQ result except that in this case if $\lambda \in [\sum_{k=1}^{j-1} \mu_k\gamma_k, \sum_{k=1}^{j} \mu_k\gamma_k)$, then all servers
in pools $k \in [j-1]$ have exactly one job and
in pool $j$ a fraction $(\lambda-\sum_{k=1}^{j-1}\mu_k \gamma_k)/\mu_j \gamma_j$ of servers have exactly one job; all remaining servers are idle. Thus, the total (scaled) number of jobs in state $x^*$ is equal to $z^*$ as defined
by~\eqref{eq:def_q*}. Thus, by Little's law, the mean response time
of jobs $\bar T_N$ under the SA-JSQ policy converges to $z^*/\lambda$ as $N \to \infty$,
which by Corollary~\ref{cor:optimality}
implies the asymptotic optimality of
the SA-JSQ policy.


\section{Stability}
\label{sec:Stability_SA-JSQ}

To show that the process $\mf Q^{(N)}$
is positive recurrent
for all $\lambda < 1$, we use an
appropriately defined Lyapunov function
and show that its drift along
any trajectory of $\mf Q^{(N)}$
is negative outside a compact subset
of the state space.
For any function $f:\mb{Z}_+^N \to \mb{R}$,
defined on the state space of 
the process $\mf Q^{(N)}$,
the drift evaluated at a state
$\mf Q \in \mb{Z}_+^N$  is 
given by
\begin{align}
\label{eq:genq_def}
    G_{\mf Q^{(N)}} f(\mf Q) = \sum_{j \in [M]}\sum_{k\in [N\gamma_j]} [& r_{k,j}^{+,N}(\mf Q)(f(\mf Q+\mf e_{k,j}^{(N)})-f(\mf Q))\nonumber \\
    & +r_{k,j}^{-,N}(\mf Q)(f(\mf Q-\mf e_{k,j}^{(N)})-f(\mf Q))],
\end{align}
where $G_{\mf Q^{(N)}}$ is the generator of $\mf Q^{(N)}$;
$\mf e_{k,j}^{(N)}$ denotes the N-dimensional 
unit vector with one
in the $(k,j)^{\textrm{th}}$ position; 
$r^{\pm,N}_{k,j}(\mf Q)$
are the transition rates from the state $\mf Q$ to the states $\mf Q \pm \mf e_{k,j}^{(N)}$. 
Intuitively, the drift, defined above, represents
the expected infinitesimal 
rate at which $f(\mf Q^{(N)}(t))$ changes when
$\mf Q^{(N)}(t)=\mf Q$.
Under the SA-JSQ
policy, for each state $\mf Q \in \Z_+^N$
and each $k \in [N\gamma_j], j \in [M]$
we have
\begin{align}
    r^{+,N}_{k,j}(\mf Q) &= \begin{cases} \frac{N\lambda}{\abs{I_{\min,j}(\mf Q)}}, &\text{if } j=j^\dag(\mf Q) \text { and } k\in I_{\min,j}(\mf Q)\label{eq:uprate_q}\\
    0, &\text{otherwise}\end{cases}\\    
    r^{-,N}_{k,j}(\mf Q) &= \mu_j \indic{Q_{k,j} > 0}, \label{eq:downrate_q}
\end{align}
where $I_{\min,j}(\mf Q)$ denotes
the set of servers with the minimum queue length in pool $j$ 
in state $\mf Q$ and
$j^\dag(\mf Q)$ denotes the fastest pool that contains 
a server with the minimum queue length in state $\mf Q$. 
The upward transition rate
$r^{+,N}_{k,j}(\mf Q)$ is obtained by multiplying the total
arrival rate with the probability that the $k^{\textrm{th}}$
server in the $j^{th}$ pool receives an arrival.
Similarly, the downward transition rate 
$r^{-,N}_{k,j}(\mf Q)$ is simply
the service rate of the $k^{\textrm{th}}$
server in the $j^{th}$ pool if the server is busy. 

{\bf Proof of Theorem~\ref{thm:SA-JSQ_Stability}}:
To prove Theorem~\ref{thm:SA-JSQ_Stability},
we compute the drift of the Lyapunov function
$\Phi:\mb{Z}_+^N \to [0,\infty)$ defined as
follows
\begin{equation}
    \Phi(\mf Q)=\sum_{j \in [M]} \sum_{k \in [N\gamma_j]} Q_{k,j}^2.
\end{equation}
%
From~\eqref{eq:genq_def} 
we have that for any $\mf Q \in \mb Z_+^N$
\begin{align}
    G_{\mf Q^{(N)}} \Phi(\mf Q) & =\sum_{j \in [M]}\sum_{k \in [N\gamma_j]} r_{k,j}^{+,N}(\mf Q)(2Q_{k,j}+1)+r_{k,j}^{-,N}(\mf Q)(-2Q_{k,j}+1) \nonumber\\
    &= 2 \sum_{j \in [M]} \sum_{k \in [N \gamma_j]} r_{k,j}^{+,N}(\mf Q) Q_{k,j}+N\lambda-2\sum_{j \in [M]}\sum_{k \in [N\gamma_j]}r^{-,N}_{k,j}(\mf Q)Q_{k,j}\nonumber \\
    &\hspace{1cm}+\sum_{j \in [M]}\mu_j B_j(\mf Q), \label{eq:phi_temp}
\end{align}
where $B_j(\mf Q)$ denotes the number of busy servers
in pool $j$ in state $\mf Q$. 
In the second equality,
we have used the facts that $\sum_{j,k}r^{+,N}_{k,j}(\mf Q)=N\lambda$ and $\sum_{k}r^{-,N}_{k,j}(\mf Q)=\mu_j B_j(\mf Q)$
which follow easily from~\eqref{eq:uprate_q} and~\eqref{eq:downrate_q}, respectively. 
We further note from~\eqref{eq:uprate_q} and~\eqref{eq:downrate_q} that $\sum_{k,j}r^{+,N}_{k,j}(\mf Q) Q_{k,j}= N \lambda Q_\min$ and $\sum_{k}r^{-,N}_{k}(\mf Q) Q_{k,j}=\mu_j \sum_{k} Q_{k,j}$, where $Q_\min$ denotes the minimum queue length in
state $\mf Q$.
Hence, from~\eqref{eq:phi_temp} we have
\begin{align}
    G_{\mf Q^{(N)}}\Phi(\mf Q) &= 2 N\lambda Q_\min - 2
    \sum_j \mu_j \sum_{k}Q_{k,j} + N\lambda + \sum_{j} \mu_j B_j(\mf Q) \label{eq:gen_temp1}\\
    & \leq -2(1-\lambda)\sum_j \mu_j \sum_{k}Q_{k,j} + N\lambda+ N.
\end{align}
In the second line, we have used the facts
that $\sum_j \mu_j \sum_{k} Q_{k,j} \geq \sum_j \mu_j \sum_{k} Q_\min = N Q_\min$ and $B_j(\mf Q) \leq N \gamma_j$. 
Hence, it follows from the above that 
if $\lambda < 1$, then the drift is strictly negative whenever 
$\sum_{j} \mu_j \sum_{k} Q_{k,j} > \frac{N(1+\lambda)}{1-\lambda},$
and is bounded by $N(1+\lambda)$ otherwise.
Thus, using the Foster-Lyapunov criterion for positive recurrence (see, e.g., Proposition D.3 of \cite{Kelly_book})
we conclude that $\mf Q^{(N)}$
is positive recurrent.\qed

\section{Uniform bounds and tightness}
\label{sec:uniform_bounds}

In this section, we first prove Theorem~\ref{thm:tail_bound}
which shows
that the stationary queue length distribution has a uniformly decaying tail for all system
sizes. 
We then use this uniform bound to establish the tightness of the sequence $(\mf x^{(N)}(\infty))_N$ in $S$.

For homogeneous systems working under the classical JSQ policy, uniform
bounds, similar to~\eqref{eq:tail_bound} can be obtained by coupling the system with another similar system
working under a `worse performing' policy such as the {\em random policy} in which the destination server for each job is chosen uniformly at random
from the set of all servers. A coupling
similar to the one described in~\cite{Turner1998} ensures that the total number of jobs
in the first system is always smaller than that in the second system.
However, in the heterogeneous setting,
it is difficult to construct a similar coupling since
the arrival of a job to a given pool
in one system does not guarantee
that the job will join the same pool
in the other system.

To overcome this difficulty, we use
a completely different approach based
on analysing the drift of an exponential
Lyapunov function. In particular,
we analyse the drift of the Lyapunov function
$
\Psi_\theta:\Z_+^N \to \R_+$ defined as
\begin{equation}
\label{eq:psi_def}
    \Psi_\theta(\mf Q) \define \sum_{j\in [M]}\sum_{k \in [N\gamma_j]} \exp(\theta Q_{k,j}), 
\end{equation}
for some $\theta > 0$. They key idea we employ
to prove Theorem~\ref{thm:tail_bound} is the following:
we show that for some positive values of $\theta$ 
the steady-state expected drift $\E[G_{\mf Q^{(N)}}\Psi_\theta(\mf Q^{(N)}(\infty))]$ of $\Psi_\theta$ is non-negative.
From this, we obtain bounds on the 
weighted sum of moment generating functions
of the stationary queue lengths of different pools, i.e.,
on $\E[\sum_{j \in [M]} \mu_j \gamma_j \exp(\theta Q^{(N)}_{k,j}(\infty))]$ for some positive $\theta$.
Using Chernoff bounds, we then obtain the bounds on 
the tail probabilities as in Theorem~\ref{thm:tail_bound}.

Note that although we would normally expect the steady-state expected drift to satisfy $$\E[G_{\mf Q^{(N)}}\Psi_\theta(\mf Q^{(N)}(\infty))] = 0,$$ proving this requires showing
$\E[\Psi_\theta(\mf Q^{(N)}(\infty))] < \infty$ which
is only true for some $\theta > 0$ if all moments of the stationary queue lengths
exist. However,
the stability condition only guarantees the 
existence of the first moment, but it does not guarantee
the existence of higher moments of the queue lengths.
Thus, to prove Theorem~\ref{thm:tail_bound}, we use
a weaker condition given by Proposition 1 of~\cite{Glynn_bounds} which states
that for any non-negative function $f:\Z_+^N \to \R_+$,
if $\sup_{\mf Q \in \Z_+^N} G_{\mf Q^{(N)}} f(\mf Q) < \infty$ then
$\E[G_{\mf Q^{(N)}} f(\mf Q^{(N)}(\infty))] \geq 0$. However, to use this result we first need to show that
the function $\Psi_\theta$ satisfies the above condition for
some positive $\theta$. This is shown in the following lemma.

\begin{lem}
For $\Psi_\theta$ as defined in \eqref{eq:psi_def} and $\theta \geq 0$, we have

\begin{equation}
\label{eq:psi_temp1}
    G_{\mf Q^{(N)}}\Psi_\theta(\mf Q) \leq (1-e^{-\theta})\sbrac{(\lambda e^\theta-1)\sum_{j \in [M]}\sum_{k \in [N \gamma_j]} \mu_j \exp(\theta Q_{k,j})+\sum_{j \in [M]} \mu_j I_j(\mf Q)},
\end{equation}
where $I_j(\mf Q)$ denotes the number of idle servers in
pool $j \in [M]$ in state $\mf Q$. In particular, for
all $\theta \in [0, -\log \lambda)$ we have
\begin{equation}
    \sup_{\mf Q \in \Z_+^N} G_{\mf Q^{(N)}}\Psi_\theta(\mf Q) < \infty,
\end{equation}
and
the steady-state drift of $\Psi_\theta$
satisfies

\begin{equation}
\label{eq:drift_bound}
    \E[G_{\mf Q^{(N)}} \Psi_\theta(\mf Q^{(N)}(\infty))] \geq 0.
\end{equation}
\end{lem}
\begin{proof}
From~\eqref{eq:genq_def} and~\eqref{eq:psi_def},
we have that for any $\mf Q \in \Z_+^N$
\begin{align}
    G_{\mf Q^{(N)}} \Psi_\theta(\mf Q)  = \sum_{j \in [M]} \sum_{k \in [N\gamma_j]} \Big[& r_{k,j}^{+,N}(\mf Q)\Big(\exp(\theta Q_{k,j})(e^\theta-1)\Big)\nonumber \\
   &+r_{k,j}^{-,N}(\mf Q)\Big(\exp(\theta Q_{k,j})(e^{-\theta}-1)\Big)\Big].
\end{align}
First, note that from~\eqref{eq:uprate_q}, we can write
\begin{align}
&\sum_{j \in [M]} \sum_{k \in [N\gamma_{j^\dag(\mf Q)}]}r_{k,j}^{+,N}(\mf Q)
\exp(\theta Q_{k,j})\nonumber \\
&=\sum_{k \in [N\gamma_{j^\dag(\mf Q)}]}\frac{ N \lambda \indic{k\in I_{\min,j^\dag(\mf Q)}(\mf Q)}\exp(\theta Q_{k,j^\dag(\mf Q)})}{\abs{I_{\min,j^\dag(\mf Q)}(\mf Q)}}, \nonumber \\
&= N \lambda \sum_{k \in [N\gamma_{j^\dag(\mf Q)}]}\frac{\indic{k\in I_{\min,j^\dag(\mf Q)}(\mf Q)}\exp(\theta Q_{\min})}{\abs{I_{\min,j^\dag(\mf Q)}(\mf Q)}},\nonumber \\
&=N\lambda\exp(\theta Q_{\min}), \label{eqn:uni_bound_int_1}
\end{align}
where $Q_{\min}$ denotes the minimum
queue length in state $\mf Q$ and the equality
on the second line follows from the fact
that for $k \in I_{\min,j^\dag(\mf Q)}(\mf Q)$ the queue
length at the $k^{\textrm{th}}$ server
of type $j^\dag(\mf Q)$ is $Q_{\min}$.
Furthermore, using~\eqref{eq:downrate_q}, we can write
\begin{align}
\sum_{j \in [M]} \sum_{k \in [N\gamma_j]}r_{k,j}^{-,N}(\mf Q)
\exp(\theta Q_{k,j})&= \sum_{j \in [M]} \sum_{k \in [N\gamma_j]} \mu_j \indic{Q_{k,j} > 0} \exp(\theta Q_{k,j}),\nonumber \\
&=\sum_{j \in [M]}\mu_j \sum_{k \in [N\gamma_j]}\Big(1-\indic{Q_{k,j}=0}\Big)\exp(\theta Q_{k,j}),\nonumber \\
&=\sum_{j \in [M]}\mu_j\sum_{k \in [N\gamma_j]}\exp(\theta Q_{k,j})-\sum_{j}\mu_j I_j(\mf Q).
\label{eqn:uni_bound_int_2}
\end{align}

Hence, using~\eqref{eqn:uni_bound_int_1} and~\eqref{eqn:uni_bound_int_2}, we can write
\begin{align}
\label{eq:psi_temp}
    G_{\mf Q^{(N)}} \Psi_\theta(\mf Q)= (e^\theta-1)\Big[ N\lambda\exp(\theta Q_\min)&-\frac{1}{e^\theta}\sum_{j\in [M]}\mu_j\sum_{k\in [N\gamma_j]}\exp(\theta Q_{k,j}) \nonumber\\
    &+\frac{1}{e^\theta}\sum_j\mu_jI_j(\mf Q)\Big]. 
\end{align}
We further note that for any $\theta >0$ we have
\begin{align*}
 \sum_j\mu_j\sum_k\exp(\theta Q_{k,j})
&\geq \sum_j\mu_j\sum_k\exp(\theta Q_\min)\\
&=N\exp(\theta Q_\min)\sum_{j}\mu_j\gamma_j\\
&=N\exp(\theta Q_\min).
\end{align*}
Using this fact in~\eqref{eq:psi_temp}, we have
that for any $\theta \geq 0$
\eq{
\begin{aligned}
    G_{\mf Q^{(N)}}\Psi_\theta(\mf Q) &\leq (1-e^{-\theta})\sbrac{\brac{\lambda e^\theta-1}\sum_{j\in [M]}\mu_j\sum_{k\in [N\gamma_j]}\exp(\theta Q_{k,j})+\sum_{j\in [M]}\mu_jI_j(\mf Q)}.
\end{aligned}
}
This proves the first statement of the lemma.
In order to prove the next statement, note that for all $\theta \in [0, -\log \lambda)$, we have $\lambda e^\theta-1<0$. Therefore, from~\eqref{eq:psi_temp1} it follows that for all
$\theta \in [0, -\log \lambda)$  we have 
\begin{align}
    G_{\mf Q^{(N)}}\Psi_\theta(\mf Q) &\leq (1-e^{-\theta})\sum_{j\in[M]}\mu_jI_j(\mf Q)\leq (1-e^{-\theta})N,
\end{align}
where in the second inequality we have used the fact that $I_j(\mf Q)\leq N \gamma_j$. Hence, 
 for all $\theta \in [0, -\log \lambda)$ we have
$
    \sup_{\mf Q \in \Z_+^N} G_{\mf Q^{(N)}}\Psi_\theta(\mf Q) < \infty
$
The last statement of the lemma now follows from the application of Proposition~1 of~\cite{Glynn_bounds}.
\end{proof}

{\bf Proof of Theorem~\ref{thm:tail_bound}}: We are
now equipped to prove Theorem~\ref{thm:tail_bound}
using the result of the lemma above. Taking expectation
of~\eqref{eq:psi_temp1} with respect to the stationary distribution and applying~\eqref{eq:drift_bound} we obtain
%
%

\begin{align}
    \brac{1-\lambda e^\theta}\E\sbrac{\sum_j\mu_j\sum_k\exp(\theta Q_{k,j}^{(N)}(\infty))}
    &\leq\E\sbrac{\sum_j\mu_j(N\gamma_j-B_j(\mf Q^{(N)}(\infty)))}\nonumber\\
    &= N(1-\lambda), \label{eqn:MGF_temp1}
    \end{align}
where $B_j(\mf Q)=N\gamma_j-I_j(\mf Q)$ denotes the number of busy servers in pool $j$ in state $\mf Q$. 
In the last inequality we have used the fact
that due to erogodicity of the process $\mf Q^{(N)}$,
the steady state rate of departure from
the system $\E[\sum_{j} \mu_j B_j(\mf Q)]$
is equal to the arrival rate $N\lambda$. Hence, from~\eqref{eqn:MGF_temp1} we 
have that for all $\theta \in [0,-\log \lambda)$

\begin{align*}
    \frac{N(1-\lambda)}{1-\lambda e^\theta} &\geq \E\sbrac{\sum_{j\in [M]}\mu_j\sum_{k\in [N\gamma_j]}\exp(\theta Q_{k,j}^{(N)}(\infty))}\\
    &= N\E\sbrac{\sum_{j\in[M]}\mu_j\gamma_j\exp(\theta Q_{k,j}^{(N)}(\infty))}\\
    &\geq N\E\sbrac{\mu_j\gamma_j\exp(\theta Q_{k,j}^{(N)}(\infty))}
\end{align*}
where the second equality follows from
the exchangeability of the stationary measure.
Thus, for each $j \in [M]$
and $\theta \in [0, -\log \lambda)$ we have

\begin{equation}
    \E\sbrac{\exp(\theta Q_{k,j}^{(N)}(\infty))} \leq \frac{1}{\mu_j\gamma_j}\frac{1-\lambda}{1-\lambda e^\theta}.
    \label{eqn:mgf_bound_1}
\end{equation}
Now, for each positive $\theta$ we have
\eq{
\mathbb{P}(Q_{k,j}^{(N)}(\infty)\geq l)=\mathbb{P}(\exp(\theta Q_{k,j}^{(N)}(\infty))\geq \exp(\theta l))\leq\frac{\E \sbrac{\exp{\big(\theta Q^{(N)}_{k,j}(\infty)}\big)}}{\exp(\theta l)}.
}
The statement of the theorem now follows by using~\eqref{eqn:mgf_bound_1} on the above
inequality.
\qed

Using Theorem~\ref{thm:tail_bound},
we now show that the sequence of stationary states
$(\mf x^{(N)}(\infty))_N$
is tight in the space $S$ under the $\ell_1$-norm.
According to Prohorov's theorem~\cite{Billingsley2013}
the tightness of this sequence will imply that the sequence
has convergent subsequences with limits in $S$. 
We shall show later that all convergent subsequences of $(\mf x^{(N)}(\infty))_N$ have the same limit,
thereby establishing the convergence of the original sequence $(\mf x^{(N)}(\infty))_N$ to the same limit.

\begin{proposition}
\label{prop:tightness}
The sequence $(\mf x^{(N)}(\infty))_N$ of stationary states
is tight in $S$ under the $\ell_1$-norm.
\end{proposition}

The necessary and sufficient criterion for tightness 
of the sequence $(\mf x^{(N)}(\infty))_N$ in $S$ under the $\ell_1$-norm 
is derived in Lemma~\ref{lem:tightness_condition}
of~\ref{sec:compact_sets}
and is given by
\eqn{
\label{eqn:l1_tight_3}
\lim_{l\to \infty} \limsup_{N\to \infty} \mathbb{P} \Big(\max_{j\in[M]}\sum_{i \geq l}x_{i,j}^{(N)}(\infty)>\epsilon\Big)=0, \ \forall \epsilon>0.
}
The proof of Proposition~\ref{prop:tightness} consists of verifying this condition using the uniform bounds derived in Theorem~\ref{thm:tail_bound}. The formal proof is given in~\ref{app:tightness}.

\section{Monotonicity}
\label{sec:monotone}

In this section, we prove the monotonicity
property stated in Theorem~\ref{thm:SA-JSQ_monotone}.
The key idea here is to couple
the arrivals and the departures
of the two systems such that if
the inequality $\mf Q^{(N)}(t^-) \leq \mf{\tilde Q}^{(N)}(t^-)$ is satisfied
just before the arrival or departure event at time $t$ then it is also
satisfied after the event has taken place.


{\bf Proof of Theorem~\ref{thm:SA-JSQ_monotone}}:
We refer to the systems corresponding to the 
processes $\mathbf{Q}^{(N)}$ and $\tilde{\mathbf{Q}}^{(N)}$
as the \textit{smaller} and \textit{larger} systems, respectively. Furthermore, in both systems,
we call the $k^{\textrm{th}}$ server
in the $j^{\textrm{th}}$ pool as the 
server with index $(k,j)$.

To couple the departures, we first generate a sequence of potential departure instants at the points of a Poisson process with rate $N$. At each potential departure instant,
a server index $(k,j)$
is chosen as follows:
First, a server type $j \in [M]$ is
chosen with probability $\mu_j\gamma_j$ (recall that $\sum_{j\in [M]} \mu_j \gamma_j=1$). 
Then, any server $k \in [N \gamma_j]$ within 
the chosen type $j$ is selected uniformly at random. 
Once the server index $(k,j)$ is chosen as described above,
the server with the index $(k,j)$ is selected 
for departure in both systems.
In each system, an actual departure occurs
from the selected server if the selected
server is busy; otherwise, no departure occurs
from the selected server 
(this is why the term {\em potential departure}
is used to describe the event).
Let $(k,j)$ denote the index of the chosen server
and $D$ denote the potential departure instant.
Assume that the inequality $Q_{k,j}^{(N)}(D^-) \leq \tilde{Q}^{(N)}_{k,j}(D^-)$ holds just before the departure.
Then,
we must have $Q_{k,j}^{(N)}(D) \leq \tilde{Q}^{(N)}_{k,j}(D)$.
This is 
because $Q_{k,j}^{(N)}(D^-) > 0$ only if $\tilde{Q}^{(N)}_{k,j}(D^-) >0$ in which case a departure occurs from both systems.

To couple the arrivals,
we generate a common Poisson arrival stream with rate $N\lambda$ for both systems. 
At each arrival instant $A$,
the job assignment decision is made
following the SA-JSQ policy
in each system independently
of the other system,
unless the pool containing
the destination server
is the same for both systems.
If $j^\dag(\mf Q^{(N)}(A^-))=j^\dag(\tilde{\mf Q}^{(N)}(A^-))=j$,
then we perform the following steps sequentially: (1) In the smaller system, the destination server is chosen
uniformly at random from the servers having 
the minimum queue length in pool $j$. Let
the index of this chosen server be $(k,j)$.
(2) We check if the server having the same index $(k,j)$
in the larger system has the minimum queue length.
If so, then this server is chosen as the destination
server in the larger system. (3) Otherwise,
the destination server in the larger system is chosen 
uniformly at random (and independently of the smaller system) 
from the servers having the minimum
queue length in pool $j$ in the larger system.
Let $(k_s,j_s)$ and $(k_l,j_l)$ be the indices of the
destination servers in the smaller and the larger systems,
respectively, at an arrival instant $A$.
Let $S=Q^{(N)}_{k_s,j_s}(A^-)$ and $L=\tilde{Q}^{(N)}_{k_l,j_l}(A^-)$.
Assume that $\mf Q^{(N)}(t) \leq \tilde{\mf Q}^{(N)}(t)$
holds for all $t < A$. Hence, $S \leq L$.
It is sufficient to show that the inequality 
\begin{equation}
    \mf Q^{(N)}(t) \leq \tilde{\mf Q}^{(N)}(t) \text{ holds at } t=A.
    \label{eq:dominance}
\end{equation}
If $(k_s,j_s)=(k_l,j_l)$, then the inequality 
trivially holds at $t=A$.
If $(k_s,j_s)\neq(k_l,j_l)$,
then we have the following possibilities:

\begin{enumerate}
    \item {\em If $j_s<j_l$}: In the larger system, the incoming arrival joins the server with index $(k_l,j_l)$ and there is no arrival to the server with index $(k_l,j_l)$ in smaller system. Hence, after the arrival we have $$Q_{k_l,j_l}^{(N)}(A)=Q_{k_l,j_l}^{(N)}(A^-)<\tilde{Q}_{k_l,j_l}^{(N)}(A^-)+1=\tilde{Q}_{k_l,j_l}^{(N)}(A).$$
    Thus, for~\eqref{eq:dominance} to be violated
    we must have $Q_{k_s,j_s}^{(N)}(A^-)=\tilde{Q}_{k_s,j_s}^{(N)}(A^-)=S$. But since $S \leq L$, we must have
    $\tilde{Q}_{k_s,j_s}^{(N)}(A^-) \leq \tilde{Q}_{k_l,j_l}^{(N)}(A^-)$ which
    contradicts with the fact that the server
    with index $(k_l,j_l)$ is the destination server
    in the larger system, i.e., it is the fastest server
    with the minimum queue length in the larger system,
    because a faster server of type $j_s (< j_l)$ has
    smaller queue length.
    Hence, \eqref{eq:dominance} must hold in this case.

    \item {\em If $j_s=j_l=j, k_s\neq k_l$}: Similarly as before, the inequality $Q^{(N)}_{k_l,j}(A) < \tilde{Q}^{(N)}_{k_l,j}(A)$ holds after the arrival. For violation of~\eqref{eq:dominance}, we therefore must have $Q_{k_s,j}^{(N)}(A^-)=\tilde{Q}_{k_s,j}^{(N)}(A^-)=S$ before the arrival. Since the arrival is assigned to the server with index $(k_l,j)$
    in the larger system, we must have $L=\tilde{Q}_{k_l,j}^{(N)}(A^-) \leq \tilde{Q}_{k_s,j}^{(N)}(A^-)=S$. Hence, we must have $S=L$.
    Also, since $S=Q^{(N)}_{k_s,j}(A^-) \leq {Q}_{k_l,j}^{(N)}(A^-)\leq \tilde{Q}_{k_l,j}^{(N)}(A^-)=L$, we must have ${Q}_{k_l,j}^{(N)}(A^-)=S=L$.
    Thus, we have $S=Q_{k_s,j}^{(N)}(A^-)=Q_{k_l,j}^{(N)}(A^-)=\tilde{Q}_{k_s,j}^{(N)}(A^-)=\tilde{Q}_{k_l,j}^{(N)}(A^-)=L$. In this case,
    our coupling rule dictates that if we have chosen the index $(k_s,j)$
    as the destination server for the smaller system, then we must
    also choose the same indexed server in the larger system as
    the destination server. This leads to a contradiction because $k_s \neq k_l$. Hence, \eqref{eq:dominance} must hold in this case.
    
     \item {\em If $j_s>j_l$}: Similar to the first two cases, the inequality $Q^{(N)}_{k_l,j_l}(A)< \tilde{Q}^{(N)}_{k_l,j_l}(A)$ holds. For violation of \eqref{eq:dominance}, we must have $$Q^{(N)}_{k_s,j_s}(A^-)= \tilde{Q}^{(N)}_{k_s,j_s}(A^-)=S$$ just before the arrival. 
     Since $L=\tilde{Q}^{(N)}_{k_l,j_l}(A^-) \leq \tilde{Q}^{(N)}_{k_s,j_s}(A^-)=S$,
     this implies $S=L$.
     Furthermore, since $S=Q^{(N)}_{k_s,j_s}(A^-) \leq {Q}_{k_l,j_l}^{(N)}(A^-)\leq \tilde{Q}_{k_l,j_l}^{(N)}(A^-)=L$,
     we must have 
     $$S=Q_{k_s,j_s}^{(N)}(A^-)=Q_{k_l,j_l}^{(N)}(A^-)=\tilde{Q}_{k_s,j_s}^{(N)}(A^-)=\tilde{Q}_{k_l,j_l}^{(N)}(A^-)=L.$$
     Hence, there must be a tie between the servers of indices $(k_s,j_s)$ and $(k_l,j_l)$ in both systems. Therefore, according to the SA-JSQ policy the incoming arrival should be assigned to the server having index $(i_l,j_l)$ in smaller system because $j_l < j_s$. This leads to the contradiction because job has been assigned to the server having index $(i_s,j_s)$ in the smaller system. Hence, \eqref{eq:dominance} must hold in this case.
\end{enumerate}
This completes the proof of the theorem. \qed



\section{Process convergence}
\label{subsec:Process_Level_Convergence_SA-JSQ}

In this section, we outline the proof of Theorem~\ref{thm:SA-JSQ_Process_level_convergence}
using the martingale approach of~\cite{Whitt2007}
and the time-scale separation technique of~\cite{Hunt1994}. 
Here, we discuss the main steps of the proof
and characterise the fluid limit process
$\mf x$.
The proof consists of the following three steps

{\bf Step 1: Martingale representation}: The first step is to express
the evolution of each component of
the process $\mf x^{(N)}$ in terms of
suitably defined martingales and a process $\mf V^{(N)}$
which evolves at a faster time scale than the components
of $\mf x^{(N)}$. In particular, the process $\mathbf{V}^{(N)}=(\mathbf{V}^{(N)}(t)=(V^{(N)}_{i,j}(t), i\geq1, j\in [M]),t\geq0)$ is defined as $V^{(N)}_{i,j}(t)=N\gamma_j-N\gamma_j x^{(N)}_{i,j}(t)$.
Thus, $V^{(N)}_{i,j}(t)$ counts the number of type $j$ servers with at most $i-1$ jobs at time $t$. It is easy to see that $\mf V^{(N)}$ is a Markov process defined on
$\mathbf{E}=(\Bar{\Z}^{\infty}_+)^M$
with transition rates 
\EQN{
\label{eqn:rate_transition_V}
\mathbf{V}^{(N)} \rightarrow
\begin{cases}
                                      \mathbf{V}^{(N)} + \mathbf{e}_{i,j},& \text{at rate $N\gamma_j\mu_j(x^{(N)}_{i,j}-x^{(N)}_{i+1,j})$} \\
                                   \mathbf{V}^{(N)} - \mathbf{e}_{i,j},& \text{at rate $N\lambda \indic{\mathbf{V}^{(N)}\in R_{i,j}}$}
\end{cases}, \ i\geq1, \ j\in[M],
}
where $R_{i,j}$ for all $i\geq1$ and $j \in [M]$ is defined as
\begin{align}
\label{eqn:R_1}
R_{i,j}=\Big\{\mathbf{v}=(v_{i,j})\in \mf E :v_{i,k}=0 \ \forall k\in [j-1],
v_{i-1,l}=0 \ &\forall l\in\cbrac{j+1,\dots,M},\nonumber\\
&0=v_{i-1,j}<v_{i,j}, \Big\}.
\end{align}
The set $R_{i,j}$ represents the set of states where the minimum queue length is at least $i$ for the pools in $\{1,\ldots,j-1\}$, exactly $i-1$
for pool $j$, and at least $i-1$ for the pools in $\cbrac{j+1,\dots,M}$. Thus, when $V^{(N)}(t) \in R_{i,j}$,
an incoming job under the SA-JSQ scheme will be assigned
to a type $j$ server with queue length $i-1$.
We can therefore 
express $x^{(N)}_{i,j}$
for each $i \geq 1$ and $j \in [M]$ as follows
\begin{multline}
\label{eqn:outline_1}
x^{(N)}_{i,j}(t)=x^{(N)}_{i,j}(0)+ \frac{\lambda}{\gamma_j} \int_0^t \mathbbm{1}_{\{\mathbf{V}^{(N)}(s)\in R_{i,j}\}}ds  - \mu_j\int_0^t (x^{(N)}_{i,j}(s)-x^{(N)}_{i+1,j}(s))ds\\
+\frac{1}{N\gamma_j}(M^{(A,N)}_{i,j}(t)-M^{(D,N)}_{i,j}(t)),
\end{multline}
%
where $M^{(A,N)}_{i,j}$ and $M^{(D,N)}_{i,j}$ are martingales corresponding to the arrivals and departures at
the component $x^{(N)}_{i,j}$. The precise definitions
of these martingales in terms
of the counting processes are given in are given in~\ref{app:martingle_rep}. It is important to note the difference
in the time-scale for the processes $\mf x^{(N)}$
and $\mf V^{(N)}$.
In a small interval $[t,t+\delta]$, the process $\mf V^{(N)}$ experiences ${O}(N\delta)$ transitions whereas the $\mf x^{(N)}$ changes only by ${O}(\delta)$. 
Hence, for large $N$ the process $\mf V^{(N)}$ reaches 
its steady-state while $\mf x^{(N)}$ remains almost constant
in this interval. This separation of the two time-scales becomes crucial in characterising the limit of the indicator
function $\indic{\mf V^{(N)}(s) \in R_{i,j}}$ appearing
in \eqref{eqn:outline_1}. 

Since the time-scales of $\mf V^{(N)}$ and $\mf x^{(N)}$ are different, they have different limits as $N \to \infty$. To treat them as a single object and characterise its limit, we define the joint process $(\mf x^{(N)}, \beta^{(N)})$ where $\beta^{(N)}$ is a random
measure defined on $[0,\infty)\times \mf E$ as $$\beta^{(N)}(A_1\times A_2)=\int_{A_1} \indic{V^{(N)}(s) \in A_2} ds.$$
for any $A_1 \in \mc{B}([0,\infty))$ and $A_2 \in \mc{B}(\mf E)$. Hence, \eqref{eqn:outline_1}
can be rewritten in terms of $\beta^{(N)}$ for $i\geq1$ and $j\in[M]$ as follows
\begin{multline}
\label{eqn:outline_3}
x^{(N)}_{i,j}(t)=x^{(N)}_{i,j}(0)+ \frac{\lambda}{\gamma_j}\beta^{(N)}([0,t]\times R_{i,j})   - \mu_j\int_0^t (x^{(N)}_{i,j}(s)-x^{(N)}_{i+1,j}(s))ds\\
+\frac{1}{N\gamma_j}(M^{(A,N)}_{i,j}(t)-M^{(D,N)}_{i,j}(t)).
\end{multline}
{\bf Step 2: Relative compactness}: The second
step in proving Theorem~\ref{thm:SA-JSQ_Process_level_convergence} consists of showing that the sequence of processes
$((\mf x^{(N)},\beta^{(N)}))_N$ is relatively compact in $D_S[0,\infty)\times\mc{L}_0$
where $\mc L_0$ is defined as the space of 
measures on $[0,\infty)\times \mf E$ satisfying
$\beta([0,t]\times \mf E)=t$ for each $t \geq 0$
and each $\beta \in \mc L_0$. We equip $\mc L_0$
with the topology of weak convergence of measures restricted to $[0,t]\times \mf E$ for each $t$. In the next lemma, we show that $((\mf x^{(N)},\beta^{(N)}))_N$ is a relatively compact sequence in $D_S[0,\infty)\times \mc{L}_0$
and characterise the limit of any convergent subsequence.

\begin{lem}
\label{lem:relative_compact_SA-JSQ}
If $\mathbf{x}^{N}(0) \Rightarrow \mathbf{x}(0)\in S$ as $N \to \infty$, then the sequence $((\mathbf{x}^{(N)},\beta^{(N)}))_{N}$ is relatively compact in $D_{S}[0,\infty) \times \mc{L}_0$ and the limit $(\mathbf{x},\beta)$ of any convergent subsequence satisfies for all $t\geq 0, i\geq1, j\in[M]$
\EQN{
\label{eqn:SA-JSQ_random_rep_limit}
x_{i,j}(t)=x_{i,j}(0)+ \frac{\lambda}{\gamma_j} \beta([0,t]\times R_{i,j})  -\mu_j\int_0^t (x_{i,j}(s)-x_{i+1,j}(s))ds.
}
\end{lem}
The proof of Lemma~\ref{lem:relative_compact_SA-JSQ} consists of verifying standard conditions
of relative compactness in $D_S[0,\infty)$
given in Proposition 3.2.4 of~\cite{kurtz_book}.
This is achieved using an approach similar to~\cite{Mukherjee2018}. The details are given in~\ref{app:relative_completness_proof}.

{\bf Step 3: Characterisation of the limit}: The final step in proving Theorem~\ref{thm:SA-JSQ_Process_level_convergence} is the characterisation the limit $\beta([0,t], R_{i,j})$
appearing in~\eqref{eqn:SA-JSQ_random_rep_limit}.
To do so, we define for any $\mf x \in S$ a Markov
process $\mf V_{\mf x}$ on $\mf E$
with transition rates
\EQN{
\label{eqn:SA-JSQrate_transition_V_x_outline}
\mathbf{V}_{\mathbf{x}} \rightarrow
\begin{cases}
                                      \mathbf{V}_{\mathbf{x}} + \mathbf{e}_{i,j},& \text{at rate $ \gamma_j \mu_j(x_{i,j}-x_{i+1,j})$} \\
                                   \mathbf{V}_{\mathbf{x}}- \mathbf{e}_{i,j},& \text{at rate $\lambda \indic{\mathbf{V}_{\mathbf{x}}\in R_{i,j}}$}
\end{cases}, \ i\geq1, \ j\in[M].
}
From Lemma~2 and Theorem~3 of~\cite{Hunt1994}, it follows that the limit $\beta([0,t]\times R_{i,j})$ 
satisfies
\eq{\beta([0,t]\times R_{i,j})=\int_0^{t}\pi_{\mathbf{x}(s)}({R_{i,j}})ds,\ i\geq1,j\in[M],}
where $\pi_{\mathbf{x}}$ is a stationary measure of the process $\mathbf{V}_{\mathbf{x}}$ and $\pi_{\mf x}$ satisfies the following for all $j \in [M]$
\EQN{
\label{eqn:SA-JSQ_pi_prop_outline}
\pi_{\mathbf{x}}(\{\mathbf{V} \in \mf E:V_{i,j}=\infty\})=1, \  \text{if} \  x_{i,j}<1.}
Hence, we can write~\eqref{eqn:SA-JSQ_random_rep_limit} in terms of $\pi_{\mathbf{x}}$ as 
\EQN{
\label{eq:evolution_x}
x_{i,j}(t)=x_{i,j}(0)+ \frac{\lambda}{\gamma_j}\int_0^{t}\pi_{\mathbf{x}(s)}({R_{i,j}})ds -\mu_j\int_0^t (x_{i,j}(s)-x_{i+1,j}(s))ds,\ i\ge1, \ j\in[M].
}
We set $p_{i-1,j}(\mf x)=\pi_{\mf x}(R_{i,j})$.
Hence, to complete the proof of Theorem~\ref{thm:SA-JSQ_Process_level_convergence},
it remains to show that $\mathbf{x}$ uniquely determines the stationary measure $\pi_{\mathbf{x}}$
and it satisfies the properties P\ref{p2}-P\ref{p6}
listed in Theorem~\ref{thm:SA-JSQ_Process_level_convergence}.

First note from~\eqref{eqn:R_1}
that the sets $R_{i,j}$
for $i \geq 1, j \in [M]$ form a partition of $\mf E$.
Since $\pi_{\mf x}$ is a probability measure on $\mf E$,
it follows that $\sum_{i\geq 1, j\in [M]}p_{i-1,j}(\mf x)=\sum_{i\geq 1, j \in [M]}\pi_{\mf{x}}(R_{i,j})=1$.
This proves P\ref{p2}.
Since $\norm{\mf x}_1 < \infty$ for each $\mf x \in S$,
there exists $l_j(\mf x)=\min\{i:x_{i+1,j}<1\}$ for all $j\in[M]$ and $\mf x \in S$. 
Observe that for each $j\in[M]$ we have $x_{i,j}=1$ for all $0\leq i \leq l_j(\mf x)$, and $x_{i,j}<1$ for all $i\geq l_j(\mf x)+1$. 
Therefore, from~\eqref{eqn:SA-JSQ_pi_prop_outline}, it follows that
\eqn{
\label{eqn:p3_proof}
\pi_{\mathbf{x}}\Big(\{V_{l_j(\mf x)+1,j}=V_{l_j(\mf x)+2,j}=\cdots=\infty:\forall j\in[M]\}\Big)=1.
}
Hence, from the definition of the set $R_{i,j}$
in~\eqref{eqn:R_1} it follows that $p_{i-1,j}(\mf x)=\pi_{\mathbf{x}}(R_{i,j})=0$ for all $i \geq l_j(\mf x)+2$ and $j \in [M]$. This proves the property P\ref{p3}.

To prove P\ref{p4}, we note 
that if $x_{i,j}=1$ for some $i \geq 1, j \in [M]$,
then $dx_{i,j}/dt \leq 0$. Using this fact
in~\eqref{eq:evolution_x} we conclude
that $\lambda \pi_{\mf x}(R_{i,j}) \leq \mu _j \gamma_j(x_{i,j}-x_{i+1,j})$ for all $(i,j)$ 
for which $x_{i,j}=1$. The definition
$l_j(\mf x)$ implies if $l_j(\mf x) > 0$
for some $\mf x$ and some $j$, then
that $x_{i,j}=1$ for all
$i \leq l_j(\mf x)$. Hence, $p_{i-1,j}(\mf x)=\pi_{\mf x}(R_{i,j})=0$ for $i \in [l_j(\mf x)-1]$.
This shows P\ref{p4}.
If $l_j(\mf x)>0$ for some $j\in[M]$, then it can be verified using property P\ref{p4} and the definition of $R_{i,j}$ that
\EQN{
\label{eqn:v_2}
\pi_{\textbf{x}}\Big(\{V_{1,j}=V_{2,j}=\cdots=V_{l_j(\mf x)-1,j}=0\}\Big)=1.
}

Now suppose $l_1(\mf x)=0$. From~\eqref{eqn:R_1}, it follows
that $p_{0,1}(\mf x)=\pi_{\mathbf{x}}(R_{1,1})=
\pi_{\mathbf{x}}(\{0<V_{1,1}\})=1$
since $\pi_{\mf x}(\{V_{1,1}=\infty\})=1$
according to~\eqref{eqn:SA-JSQ_pi_prop_outline}. This proves the P\ref{p5}.
The proof of P\ref{p6} follows similarly using~\eqref{eqn:R_1},~\eqref{eqn:SA-JSQ_pi_prop_outline}, and P\ref{p3}. 

Hence, the only part left 
to prove Theorem~\ref{thm:SA-JSQ_Process_level_convergence}
is to show that $\pi_{\mf x}$ is uniquely determined by $\mf x$ for all $\mf x \in S$. To show this, it is sufficient
to prove that the stationary distribution
of $(V_{l_j(\mf x),j}:j\in[M])$
is uniquely determined by $\mf x$
because the stationary distribution of all other components of
$\mf V_{\mf x}$ has already been uniquely characterised by~\eqref{eqn:p3_proof} and~\eqref{eqn:v_2}. 
The transition rates of the individual components of the chain $(V_{l_1(\mf x),1},V_{l_2(\mf x),2},\dots,V_{l_M(\mf x),M})$ are given by

\EQN{
\label{eqn:tran_rates_m_dim_markov}
V_{l_k(\mf x),k} \rightarrow
\begin{cases}
                                      V_{l_k(\mf x),k} + 1,& \text{at rate $\gamma_k \mu_k(x_{l_k(\mf x),k}-x_{l_k(\mf x)+1,k})$} \\
                                  V_{l_k(\mf x),k}- 1,& \text{at rate $\lambda \indic{\mf V \in R_{l_k(\mf x),k}}$}
\end{cases}, \forall k\in[M].
} 
Note that the Markov chain given by~\eqref{eqn:tran_rates_m_dim_markov} is defined on $\bar{\mb{Z}}_+^M$ and has $2^M$ communicating classes since each component
of the chain can be either finite or infinite.
To show the uniqueness $\pi_{\mf x}$ we need to show that $\pi_{\mf x}$ is concentrated only on a single communicating class among these $2^M$ classes. This is equivalent to showing $\pi_{\mathbf{x}}(V_{l_j(\mf x),j}=\infty)=0$ or $1$ for each $j \in [M]$. 
To show the above, we use the result of the next lemma which characterises
the stationary distribution of a finite dimensional Markov chain whose transition rates have a form similar to~\eqref{eqn:tran_rates_m_dim_markov}.

\begin{lem}
\label{lem:M-dim-Markov}
For $K \in \mb{N}$, let $\mf U=(\mf U(t)=(U_j(t), j\in [K])\in \mb{Z}_+^K:t\geq0)$ be a Markov chain with transition rates 
\EQN{
\label{eqn:k-dim-chain}
\mf U \rightarrow
\begin{cases}
                                      \mf U + \mf e_i,& \text{at rate $\nu_i$} \\
                                  \mf U - \mf e_i,& \text{at rate $\lambda \indic{0=U_{1}=\dots=U_{i-1}< U_{i}}$}
\end{cases}, \forall i\in[K],
}
where $\mf e_i$ denotes the $K$-dimensional unit vector with entry $1$ at $i^{\textrm{th}}$ position. The Markov chain $\mf U$ is positive recurrent if and only if $\sum_{i\in[K]}\frac{\nu_i}{\lambda}<1$. Furthermore, if $\pi$ denotes the stationary distribution of the chain $\mf U$, then we have 
\eqn{
\label{eqn:stat_prop_finite_M}
\pi\cbrac{0=U_{1}=\dots=U_{i-1}< U_{i}}=\frac{\nu_i}{\lambda}, \ \forall i\in[K].
}
\end{lem}
 Using Lemma~\ref{lem:M-dim-Markov}, we show that $\pi_{\mathbf{x}}(V_{l_j(\mf x),j}=\infty)=0$ or $1$ for each $j\in[M]$ in~\ref{app:Unique_charac}.
 


\section{Fixed Point Characterisation}
\label{sec:Global_stability_Limit_interchange}
In this section, we prove  Theorem~\ref{thm:SA-JSQ_global_stability_limit_interchange} which characterises the fixed point $\mf x^*$ of the fluid limit $\mf x$ and shows that the fixed point is globally attractive. For the proof of uniqueness of the fixed point, we explicitly use properties P\ref{p2}-P\ref{p6} listed in Theorem~\ref{thm:SA-JSQ_Process_level_convergence}. Moreover, to prove global stability we use the monotonicity  of the process $\mf x^{(N)}$ proved in Theorem~\ref{thm:SA-JSQ_monotone}.

{\bf Proof of Theorem~\ref{thm:SA-JSQ_global_stability_limit_interchange}.(\ref{fixed_point})}:
From~\eqref{eqn:SA-JSQ_fluid_process} it follows that for $\mathbf{x}^* \in S$ to be a fixed point of the fluid limit $\mf x$, we must have
\EQN{
\label{eqn:SA-JSQ_differential_form_1}
 \frac{\lambda}{\gamma_j}p_{i-1,j}(\mf x^*) =\mu_j(x_{i,j}^*-x_{i+1,j}^*),\ i\ge1, \ j\in[M].
}
Summing~\eqref{eqn:SA-JSQ_differential_form_1} over all $i\geq1$ and for all $j\in[M]$, we get 
$\lambda \sum_{i\geq1} \sum_{j\in [M]}p_{i-1,j}(\mf x^*)=\sum_{j\in [M]}\mu_j\gamma_jx^*_{1,j}$.
Using P\ref{p2}, this implies that 
\eqn{
\label{eqn:dep_arri_rate}
\lambda =\sum_{j\in [M]}\mu_j\gamma_jx^*_{1,j}.}
Thus, $x_{1,j}^*=1$ for all $j\in[M]$ is not possible because the stability condition requires $\lambda<1$. Hence, we must have $x_{1,j}^* < 1$
for at least one $j \in [M]$.
In the following, we consider different cases
based on the interval in which $\lambda$ belongs

    {\em If $0<\lambda<\mu_1\gamma_1$}: For $\lambda\in (0,\mu_1\gamma_1)$, we show that $x_{1,1}^*=\lambda/\mu_1\gamma_1$ and $x_{i,j}^*=0$ for all $(i,j)\neq (1,1)$. Suppose $x_{1,1}^*<1$, this means that $l_1(\mf x^*)=0$. Therefore, from property P\ref{p5}, we have $p_{0,1}(\mf x^*)=1$. Hence, summing~\eqref{eqn:SA-JSQ_differential_form_1} over all $i\geq1$ and for $j=1$, we get $x_{1,1}^*=\lambda/\mu_1\gamma_1$. Similarly, summing~\eqref{eqn:SA-JSQ_differential_form_1} for all $i\geq m$ and for $j=1$, we get $x^*_{m,1}=0$ for any $m \geq 2$. By similar line of arguments as above, we can easily verify that $x^*_{i,j}=0$ for all $i\geq1$ and for all $j\in\cbrac{2,\dots,M}$. Now, suppose $x_{1,1}^*=1$. Then from~\eqref{eqn:dep_arri_rate}, with $x_{1,1}^*=1$ implies that $\sum_{j=2}^M\mu_j\gamma_jx_{1,j}^*=\lambda-\mu_1\gamma_1<0$, which leads to a contradiction as $\mf x^*\in S$.
    
    {\em If $\sum_{i=1}^{j-1}\mu_i\gamma_i\leq\lambda<\sum_{i=1}^j\mu_i\gamma_i$, for $j\in \cbrac{2,\dots,M}$}: For this case we show that $x_{1,k}^*=1$ for all $k\in [j-1]$, $x_{1,j}^*=(\lambda-\sum_{i=1}^{j-1}\mu_i\gamma_i)/\mu_j\gamma_j$, $x_{1,k}=0$ for all $k\geq j+1$, and $x_{l,k}^*=0$ for all $k\in[M]$ and for all $l\geq2$. First, we use induction to prove that $x_{1,k}^*=1$ for all $k\in [j-1]$. Suppose $x_{1,1}^*<1$. This means that $l_1(\mf x^*)=0$. Therefore, using P\ref{p5} and summing~\eqref{eqn:SA-JSQ_differential_form_1} for all $i\geq1$ and $j=1$, we get
  $x_{1,1}^*=\lambda/\mu_1\gamma_1\geq1$,
 which contradicts the assumption that $x_{1,1}^*<1$. Therefore, $x_{1,1}^*=1$ and $l_1(\mf x)\geq1$. This proves base case for the induction. Now assume $x_{1,k}^*=1$ for all $k\in[j-2]$, which implies that $l_k(\mf x^*)\geq1$ for all $k\in[j-2]$. Using the assumption that $x_{1,k}^*=1$ for all $k\in[j-2]$, we show that $x_{1,j-1}^*=1$. Suppose $x_{1,j-1}^*<1$, which implies that $l_{j-1}(\mf x^*)=0$. Hence, using property P\ref{p6} and using~\eqref{eqn:SA-JSQ_differential_form_1} we get $\frac{\lambda}{\gamma_{j-1}\mu_{j-1}}\sum_{i\geq1}p_{i-1,j-1}(\mf x^*) =x_{1,j-1}=\frac{\lambda}{\gamma_{j-1}\mu_{j-1}}\geq1$,
    which contradicts the assumption that $x_{1,j-1}^*<1$. Therefore, we must have $x_{1,j-1}^*=1$ and $l_{j-1}(\mf x^*)\geq1$. Next, we prove that $x_{1,j}^*=(\lambda-\sum_{i=1}^{j-1}\mu_i\gamma_i)/\mu_j\gamma_j$. Suppose $l_{j}(\mf x^*)\geq1$. This implies that $x_{1,j}^*=1$. Therefore, using~\eqref{eqn:dep_arri_rate}, we have
    $\sum_{i=j+1}^M \gamma_i\mu_i x_{1,i}^*=\lambda-\sum_{i=1}^{j}\gamma_i \mu_i <0$,
    which is not possible as $\mf x^*\in S$. Hence, we have $l_{j}(\mf x^*)=0$. So far we have proved that $l_{k}(\mf x^*)\geq1$ for all $k\in[j-1]$, $l_j(\mf x^*)=0$ and $l_{k}(\mf x^*)\geq0$ for all $k\geq j+1$. Therefore, using property P\ref{p6} and equation~\eqref{eqn:SA-JSQ_differential_form_1}, we can easily get $x_{1,k}^*=0$ for all $k\geq j+1$. Now using~\eqref{eqn:dep_arri_rate}, we obtain $x_{1,j}^*=(\lambda-\sum_{i=1}^{j-1}\mu_i\gamma_i)\mu_j\gamma_j$.
    Similarly, using property P\ref{p6} and~\eqref{eqn:SA-JSQ_differential_form_1}, we can easily verify that $x_{l,k}^*=0$ for all $k\in[M]$ and for all $l\geq2$.
\qed


The proofs of global stability and limit interchange are given in~\ref{app:global_stability}. 

\section{Optimality of the Resource Pooled System}
 \label{sec:Rp_optimality}
 
In this section, we prove Theorem~\ref{thm:RP_optimality}.
We use the following lemma which establishes that
when both $\mc M'_N$ and $\mc M_N$ have the same
total number of jobs $i$, the rate of departure 
$q^{Z^{(N)}}(i,i-1)$ in $\mc M'_N$ is higher than
the rate of departure in $\mc M_N$.
Note that the departure rate in $\mc M_N$
at any time $t \geq 0$ under a policy $\Pi$ is $\sum_{j \in [M]}\mu_j X_{1,j}^{(N,\Pi)}$.

\begin{lem}
\label{lem:back_rate_thm1}
For any stationary policy $\Pi$,
if $Z^{(N)}(t)=R^{(N,\Pi)}(t)=i$ for some $t\geq0$, then

\EQN{
\label{eqn:back_rate}
\sum_{j\in[M]}\mu_j X^{(N,\Pi)}_{1,j}(t)\leq q^{Z^{(N)}}(i,i-1)=\sum_{j=1}^M \mu_j\brac{\brac{i-\sum_{i=1}^{j-1}N\gamma_i}_+ \land N\gamma_j}.
}
 \end{lem}
{\bf Proof of Theorem~\ref{thm:RP_optimality}}: 
We construct a coupling between the processes $Z^{(N)}$ and
$\mf X^{(N,\Pi)}$ such that if $Z^{(N)}(0) \leq R^{(N,\Pi)}(0)$
then $Z^{(N)}(t) \leq R^{(N,\Pi)}(t)$ for all $t \geq 0$.
Let the current instant be $t$ and assume that $Z^{(N)}(t) \leq R^{(N,\Pi)}(t)$. We shall describe a way of generating the next event 
and the time for the next event $s$ such that 
$Z^{(N)}(s) \leq R^{(N,\Pi)}(s)$ is maintained right after the event has taken place.

We generate the time until the next potential arrival for both systems
as an exponentially distributed random variable with mean $N\lambda$. Hence, for both systems, arrivals occur at the same instants.

If $Z^{(N)}(t) < R^{(N,\Pi)}(t)$ then the time until the next potential departure is generated independently for each system as exponential random variables with means $q^{Z^{(N)}}(Z^{(N)}(t),Z^{(N)}(t)-1)$ and $\sum_{j} \mu_j X_{1,j}^{(N,\Pi)}(t)$ for systems $\mc M'_N$ and $\mc M_N$, respectively. If $Z^{(N)}(t) = R^{(N,\Pi)}(t) =i$, then
we generate the time $D$ until the next potential departure for $\mc M_N$
as an exponential random variable with mean $\sum_{j} \mu_j X_{1,j}^{(N,\Pi)}(t)$. We also generate another independent
exponential random variable $C$ with mean $q^{Z^{(N)}}(i,i-1)-\sum_{j} \mu_j X_{1,j}^{(N,\Pi)}(t) \geq 0$.
Note that we can do so by Lemma~\ref{lem:back_rate_thm1}.
Now we generate the time until the next potential departure for $\mc M_N'$
as $D'=\min(D,C)$. Therefore, $D' \leq D$
and $D'$ is exponentially distributed with mean $q^{Z^{(N)}}(i,i-1)$.

Once all the potential arrival and departures are generated as described above, the next {\em true event} is taken to be the potential event which occurs the earliest; all other potential events are discarded (which is possible because of the memoryless property of the exponential distribution).
Due to the construction above, it is clear that $Z^{(N)}(s) \leq R^{(N,\Pi)}(s)$ is maintained right after the next true event.
This completes the proof. \qed

\section{Numerical Studies}
\label{sec:Numerical_Studies}
\begin{figure}
\centering
\subfigure[%
Mean response time as a function of normalized arrival rate $\lambda$ with $N=1000$ servers. We set $\mu_1=4\mu_2=20/8$,
$\gamma_1=1-\gamma_2=1/5$, and $d_1=d_2=2$. %
\label{Fig:diff_schemes}]{%
\includegraphics[width=6.5cm]{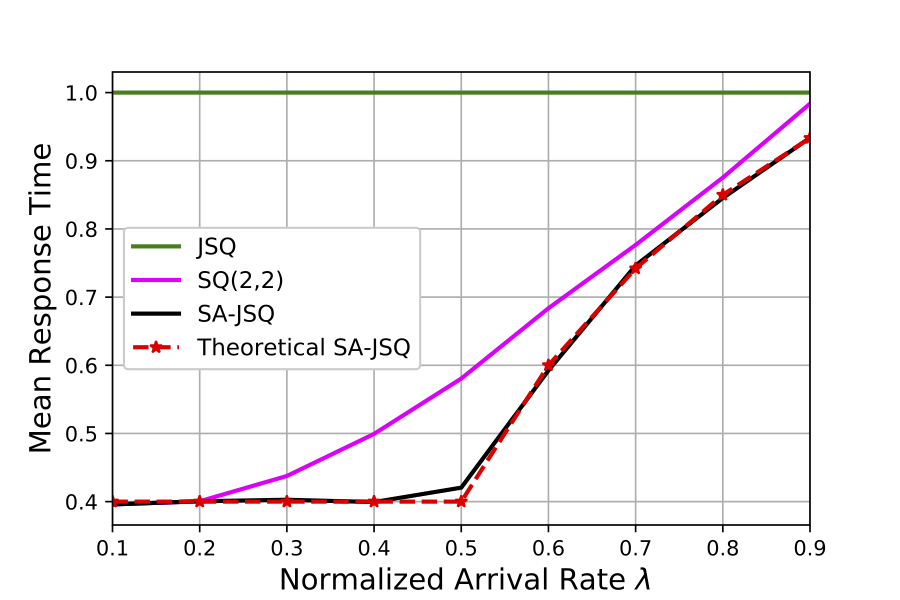}
}
\hfill
\centering
\subfigure[%
$d(\mf x^{(N)}(\infty),\mf x^*)=\sum_{i,j}|x_{i,j}^{(N)}(\infty)-x_{i,j}^*|$ as a function of system size $N$. We set $\mu_1=2\mu_2=4/3$,
$\gamma_1=1-\gamma_2=1/2$. %
\label{Fig:conv}]{%
\includegraphics[width=6.5cm]{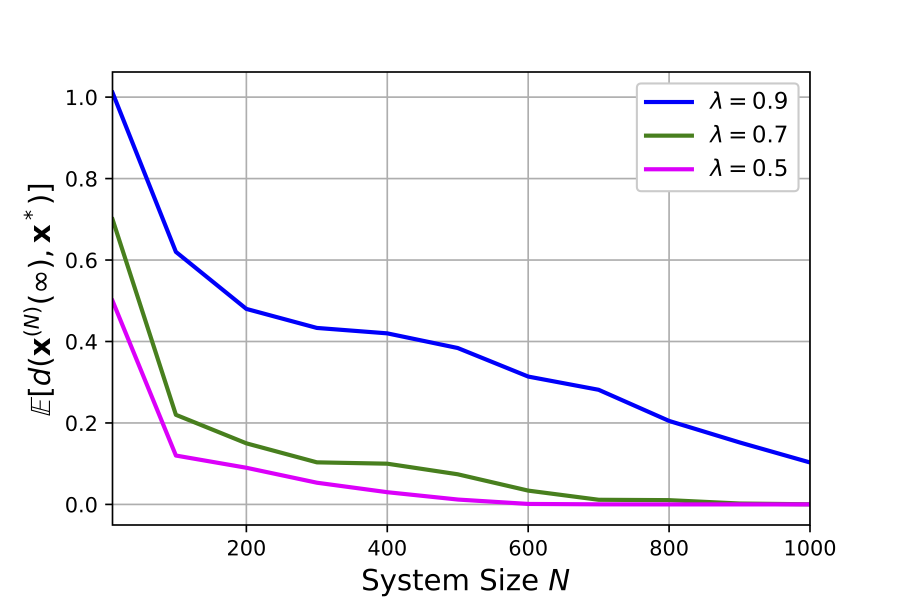}
}
\caption{Simulation plots}
\label{Fig:Cost}
\end{figure}

In this section, we present simulation results for different load balancing schemes. 
 For all simulations, we have assumed $M=2$ and taken the number of arrivals to be $3 \times10^5$.
In Figure~\ref{Fig:diff_schemes}, we have plotted the mean response time of jobs for different schemes as a function of the normalised arrival rate $\lambda$. For performance comparison we also simulated a scheme proposed in~\cite{Arpan2016} and referred to as the SQ$(d_1,d_2)$ scheme. For the SQ$(d_1,d_2)$ scheme, upon job arrival $d_j$ servers of type $j$ are sampled uniformly at random from the set of $N\gamma_j$ servers for $j\in[2]$. The job is then sent to the server with the minimum queue length among the sampled servers. Ties within different types servers are broken by selecting the type with the maximum rate. 
 We see that with SA-JSQ we obtain upto $60\%$ reduction in average response time of jobs compared to classical JSQ. As expected, the performance of SQ$(2,2)$ lies in between classical JSQ and SA-JSQ. To investigate the convergence rate to the fixed point of the fluid limit, in Figure~\ref{Fig:conv}, we have plotted the distance $d(\mf x^{(N)}(\infty),\mf x^*)=\sum_{i,j}|x_{i,j}^{(N)}(\infty)-x_{i,j}^*|$ as a function of $N$ for $\lambda\in \cbrac{0.5,0.7,0.9}$. We note that for large values of $\lambda$ the distance is higher than that for smaller values of $\lambda$.

\section{Conclusion and Future Works}
\label{sec:Conclusion}
In this paper, we have investigated speed-aware JSQ-type load balancing schemes for heterogeneous systems. We obtained a lower bound on the mean response time of jobs under any load balancing scheme by comparing the system with an appropriate resource pooled system. We showed that the  lower bound is achieved by the SA-JSQ scheme in the fluid limit, thereby establishing the asymptotic optimality of SA-JSQ. Moreover, in establishing the fluid limit of SA-JSQ, we have proved uniform bounds on the stationary measures of queue lengths which are required to prove tightness. Using coupling, we have also shown that any stationary load balancing scheme in a system of parallel queues stochastically dominates the JFFS scheme in the resource pooled system.
There are many interesting avenues for future work.
Characterising the performance of the SA-JSQ scheme in the Halfin-Whitt regime remains as an open problem. It is also interesting to analytically characterise the distance between $\mf x^{(N)}(\infty)$ and the fixed point $\mf x^*$ as a function of the system size $N$.

\bibliographystyle{IEEEtran}
\bibliography{Report1}

\appendix
\section{Resource Pooled System under JFFS: Proofs of Proposition~\ref{thm:stability_RP} and Proposition~\ref{thm:mf_rp}}
\label{sec:JFFS_Stability}

In this section, we analyse the resource pooled system $\mc M_N'$
as described in Section~\ref{Sec:RP} under the JFFS policy
and prove the results stated in Propositions~\ref{thm:stability_RP}
and~\ref{thm:mf_rp}.

\subsection{Proof of Proposition~\ref{thm:stability_RP}}

To prove the proposition, we construct a coupling between the 
resource pooled system $\mc M_N'$
and an M/M/1/N system with Poisson arrival rate $N\lambda$
and $N$ identical serves, each with rate $1$.
Let $Y^{(N)}(t)$ denote the number of jobs in the
M/M/1/N system 
at time $t\geq 0$. 
The transition rates $q^{Y^{(N)}}(k,l)$,
for $k,l \in  \mb{Z}_+$,
of the chain $Y^{(N)}=(Y^{(N)}(t), t\geq 0)$ are given
by
\begin{equation}
    q^{Y^{(N)}}(k,l)=\begin{cases}
                        N\lambda, &\text{ if } l=k+1,\\
                        k \wedge N, &\text{ if } l=k-1,\\
                        0, &\text{ otherwise}.
    \end{cases}
\end{equation}
From the standard results on
M/M/1/N queues, it follows that the process 
$Y^{(N)}=(Y^{(N)}(t), t\geq 0)$ is positive recurrent for $\lambda< 1$.
Furthermore, if $Y^{(N)}(\infty)$ denotes
the stationary number of jobs in the system,
then from the standard results on M/M/1/N queues
we have
\begin{equation}
\frac{\E\sbrac{Y^{(N)}(\infty)}}{N}=\lambda+ \frac{\lambda}{1-\lambda}\frac{\mathbb{P}\sbrac{ Y^{(N)}(\infty)\geq N}}{N} \leq \lambda +\frac{\lambda}{1-\lambda}.
\end{equation}
Hence, to prove the proposition it suffices to construct a 
coupling between $Z^{(N)}$ and $Y^{(N)}$ 
such that $Z^{(N)}(0) \leq Y^{(N)}(0)$ 
implies 
$Z^{(N)}(t) \leq Y^{(N)}(t)$ for all $t \geq 0$. 
First, note that for each Markov chain, the transition 
rate out of any state is bounded above by the constant
$B=N(\lambda+1)$.
Hence, we can generate 
both $Z^{(N)}$ and $Y^{(N)}$ processes by constructing the 
corresponding uniformized
discrete-time Markov chains $\tilde{Z}^{(N)}=(\tilde{Z}^{(N)}(m), m \in \mb{Z}_+)$
and $\tilde{Y}^{(N)}=(\tilde{Y}^{(N)}(m), m \in \mb{Z}_+)$.
The one-step transition probabilities
of these two chains from state $k$
to state $l$ are respectively given by
\begin{align*}
    p^{\tilde{Z}^{(N)}}(k,l)&=\begin{cases}
                                     \frac{q^{Z^{(N)}}(k,l)}{B}, &\text{ for } l \neq k\\
                                     1-\sum_{l'\neq k} p^{\tilde{Z}^{(N)}}(k,l'), &\text{ for } l=k,
    \end{cases}\\
    p^{\tilde{Y}^{(N)}}(k,l)&=\begin{cases}
                                     \frac{q^{ Y^{(N)}}(k,l)}{B}, &\text{ for } j \neq i\\
                                     1-\sum_{l'\neq k} p^{\tilde{Y}^{(N)}}(k,l'), &\text{ for } l=k.
    \end{cases}
\end{align*}
where $q^{Z^{(N)}}$ and  $q^{Y^{(N)}}$ denote the transition rates for $Z^{(N)}$ and $Y^{(N)}$, respectively.
To construct the continuous-time sample paths of the original
chains $Z^{(N)}$ and $Y^{(N)}$ on the same probability space,
we generate a common Poisson process with rate $B$ and embed
the time-steps of both $\tilde Z^{(N)}$ and $\tilde Y^{(N)}$
into the points of the Poisson
process.

It is easy to see from the transition rates
that $p^{\tilde{Z}^{(N)}}(k,k-1) \geq p^{\tilde{Y}^{(N)}}(k,k-1)$
for all $k \in \mb{Z}_+$. 
We let the chains $\tilde{Z}^{(N)}$ and $\tilde{Y}^{(N)}$ evolve 
independently of each other except at instants when they become equal.
If $\tilde{Z}^{(N)}(m)=\tilde{Y}^{(N)}(m)=k \in \mb{Z}_+$ for some time step $m\in \mb{Z}_+$, we first construct $\tilde{Y}^{(N)}(m+1)$ according to the transition probabilities $p^{\tilde{Y}^{(N)}}$. Then we generate $\tilde{Z}^{(N)}(m+1)$ as follows:

\eq{
    \tilde{Z}^{(N)}(m+1)=\begin{cases}
                                     \tilde{Y}^{(N)}(m+1), &\text{ if } \tilde{Y}^{(N)}(m+1)\in \{i+1,i-1\}\\
                                  \tilde{Y}^{(N)}(m)-\theta, &\text{ otherwise },   
    \end{cases}
}
where $\theta \in \{0,1\}$ is a Bernoulli random variable with $$\mathbb{P}[\theta=1\vert \tilde{Z}^{(N)}(m)=\tilde{Y}^{(N)}(m)=k]=\frac{p^{\tilde{Z}^{(N)}}(k,k-1) - p^{\tilde{Y}^{(N)}}(k,k-1)}{p^{\tilde{Y}^{(N)}}(k,k)} \geq 0.$$
Clearly, under the coupling described above $\tilde{Z}^{(N)}(m) \leq \tilde{Y}^{(N)}(m)$ for all $m \in \mb{Z}_+$ if ${Z}^{(N)}(0)=\tilde{Z}^{(N)}(0) \leq \tilde{Y}^{(N)}(0)=Y^{(N)}(0)$. Hence, we have ${Z}^{(N)}(t) \leq Y^{(N)}(t)$ for all $t \geq 0$. \qed

\subsection{Proof of Proposition~\ref{thm:mf_rp}}

To prove Proposition~\ref{thm:mf_rp},
we study the limit of the process
$z^{(N)}=\brac{z^{(N)}(t),t \geq 0}$ defined as 
\EQ{z^{(N)} (t)=\frac{Z^{(N)} (t)}{N} , \ \ t\geq 0.}
Thus, $z^{(N)}(t)$ denotes the scaled number of jobs in $\mc M'_N$ under the JFFS scheme at time $t$. 
We first characterise the limit of the sequence
of processes $(z^{(N)})_N$ in the lemma below.

\begin{lem}
\label{thm:RP_Main_thm}
If ${z}^{(N)}(0) \Rightarrow {z}(0) \in \R$ as $N \rightarrow \infty$,
then ${z}^{(N)} \Rightarrow {z}$ as $N \rightarrow \infty$, where ${z}=\brac{z(t),t \geq 0}$ is the unique deterministic process satisfying the following integral equation
\begin{equation}
\label{eqn:RP_fluid_process}
{z}(t)={z}(0) + \lambda t
-\int_0^t
\sum_{j=1}^M \mu_j\brac{\brac{{z}(s)-\sum_{i=1}^{j-1}\gamma_i}_+ \land \gamma_j}\;ds.   
\end{equation}
Furthermore, the process $z$
satisfying~\eqref{eqn:RP_fluid_process} has a unique fixed point
$z^*$ given by
\begin{equation}
\label{eqn:RP_fixed_point}
 {z}^* =  \max_{j\in [M]}\brac{
 \sum_{i=1}^{j-1}\gamma_i+ \frac{\lambda - \sum_{i=1}^{j-1}\mu_i\gamma_i}{\mu_j}}.   
\end{equation}
\end{lem}

\begin{proof}
Let $(A(t): t\geq 0)$ and $(D(t):t\geq 0)$ 
be independent unit-rate Poisson processes. 
We can express the evolution of the total number $Z^{(N)}(t)$
of jobs in the system $\mc M_N'$ as follows:
\begin{equation}
\label{eqn:app_1}
Z^{(N)} (t) = Z^{(N)} (0) + A(N\lambda t) -D(\int_0^t T(Z^{(N)}(s))\;ds), \ \ t\geq 0,    
\end{equation}
where $Z^{(N)}(0)=N z^{(N)}(0)$ and
$T: \mb R_{+} \to \mb R_+$ is defined as
\eq{T(Z)=\sum_{j=1}^M \mu_j\brac{\brac{{Z}-\sum_{i=1}^{j-1}N\gamma_i}_+ \land N\gamma_j}.}  
Hence, for $Z\in \mb{Z}_+$, $T(Z)$ represents
the total rate at which jobs depart the system in state $Z$.
Dividing~\eqref{eqn:app_1} by $N$, we have
\begin{equation}
z^{(N)} (t) = z^{(N)} (0) + M_A^{(N)}(t) -M_D^{(N)}(t) + \lambda t - \int_0^t T{'}(z^{(N)}(s)) \;ds, \ \forall t\geq0,    
\label{eq:zn_evolution}
\end{equation}
where 
\begin{align*}
 M_A^{(N)}(t)&=\frac{A(N\lambda t) -N\lambda t}{N},\\  
 M_D^{(N)}(t)&=\frac{D(\int_0^t T(Z^{(N)}(s))\;ds) -\int_0^t T(Z^{(N)}(s))\;ds}{N},\\
 T'(z)&=\frac{T(Z)}{N}=\sum_{j=1}^M \mu_j(({z}-\sum_{i=1}^{j-1}\gamma_i)_+ \land \gamma_j).
\end{align*}
Using Lemma 3.2 from~\cite{Pang2007}, it can be easily verified that $M_A^{(N)}$ and $M_D^{(N)}$ are square-integrable martingales with respect to the filtration $\mathbf{F}^{(N)}=(\mathcal{F}_{N,t},t\geq 0)$, where
\EQ{
\begin{aligned}
\mathcal{F}_{N,t} = &\sigma \Big[Z^{(N)} (0),A(N\lambda s),D\Big(\int_0^s T(Z^{(N)}(u))du\Big), 0\leq s \leq t \Big].
\end{aligned}
}
Moreover, quadratic variation process for $M_A^{(N)}$ and $M_D^{(N)}$ are given by
\EQ{
\begin{aligned}
[M^{(N)}_A](t)=& \frac{A(N\lambda t)}{N^2}, \ t\geq0\\
[M^{(N)}_D](t)=& \frac{1}{N^2}D(\int_0^t T(Z^{(N)}(s))\;ds), \ t\geq0.\\
\end{aligned}
}
Using the martingale functional central limit theorem~\cite{Whitt2007}, we have 
\EQ{M_A^{(N)} \Rightarrow 0, \ \ M_D^{(N)} \Rightarrow 0 \ \text{as} \ N \rightarrow \infty.}

From~\eqref{eq:zn_evolution} and~\eqref{eqn:RP_fluid_process},
we see that both $z^{(N)}$ and $z$
can be expressed as $z^{(N)}=f(z^{(N)}(0),M_A^{(N)}-M_D^{(N)})$
and $z=f(z(0),0)$, where $f :\R \ \times \mathcal{D}[0,\infty) \  \rightarrow \mathcal{D}[0,\infty)$ is defined as the mapping that takes $(b,y)$
to $x$ determined by the following integral equation
\EQ{
x(t)= b + y(t)+\lambda t  - \int_0^t T'(x(s))ds.
}
Hence, if $f$ is well-defined and continuous,
then the continuous mapping theorem 
proves the first statement of the lemma.
From Theorem 4.1 of~\cite{Pang2007}), it follows that
to show that $f$ is well-defined and continuous, 
it is sufficient to show that $T'$ is Lipschitz continuous.
Since, $\mu_j(({z}-\sum_{i=1}^{j-1}\gamma_i)_+ \land \gamma_j)$ is Lipschitz with constant $\mu_j$ for all $j\in [M]$. Therefore, the sum of Lipschitz functions is again Lipschitz with constant $(\mu_1\vee \mu_2 \vee \dots \vee \mu_M)$.
This establishes the first statement of the lemma.

To prove the second statement of the lemma,
note that we can can express~\eqref{eqn:RP_fluid_process} 
in its differential form as follows
\eq{\frac{d}{dt}z(t)=\lambda -T{'}(z(t)).}
Hence, any fixed point $z^*$ of the process $z$
must satisfy the equation
$\lambda -T{'}(z^*)=0$. 
Since $T'$ is a piecewise-linear map,
it is easy to solve the above equation in closed form
and find the unique solution to be~\eqref{eqn:RP_fixed_point}.
\end{proof}

\begin{lem}
\label{lem:RP_global_stability}
Let $z(u,t)$ denote the solution to~\eqref{eqn:RP_fluid_process} 
for $z(0)=u$. 
Then for any $u\in \R$, $z(u,t)\rightarrow z^*$ as $t \rightarrow \infty$.
Furthermore, the sequence $(z^{(N)}(\infty))_N$ of stationary states
converges weakly to $z^*$ as $N \to \infty$.
\end{lem}

\begin{proof}
We prove 
the first statement of the lemma by considering 
the following two cases (i) $u \geq z^*$ and (ii) $u < z^*$. 
We only provide the proof for the first case as the proof of 
the second case is similar.

If $u \geq z^*$ then $z(u,t)\geq z^*$ for all $t\geq0$. 
To see this, assume on the contrary that $z(u,t) < z^*$ for some $t > 0$.
Since $z(0,u)=u \geq z^*$ and $z(u,\cdot)$ is continuous,
there must exist a $t_1 \in (0,t)$ such that $z(u, t_1)=z^*$. 
But since $z^*$ is a fixed point, this implies that $z(t)=z^*$ for all $t\geq t_1$ which leads to a contradiction. 

Now consider the Lypunov function $V:\R_+\to \R_+$
defined as
$V(z(u,t))= z(u,t) - z^* \geq 0$. 
Using~\eqref{eqn:RP_fluid_process} we have
\EQN{
\label{eqn:app_2}
\frac{d}{dt}V(z(u,t))=\lambda - T'(z(u,t)),
}
where $T'(z)=\sum_{j=1}^M \mu_j(({z}-\sum_{i=1}^{j-1}\gamma_i)_+ \land \gamma_j)$.
Since the fixed point $z^*$ solves $\lambda - T'(z^*)=0$,
we have
\EQ{
\frac{d}{dt}V(z(u,t))=-(T'(z(u,t))-T'(z^*)).
}
We further note that $T'$ is a strictly increasing function
of its argument.
Hence, $\dot V(z(u,t)) < 0$ when $z(u,t) > z^*$.
This implies that $V(z(u,t)) \to 0$ as $t \to \infty$,
thereby proving the first part of the lemma.

To prove the second part of the lemma, we first show that the sequence $\brac{z^{(N)}(\infty)}_{N}$ is tight. 
By the application of Markov inequality and
the bound in~\eqref{eqn:sup_bound}, we have
\eq{
\sup_{N\geq1}\mathbb{P}(z^{(N)}(\infty)>a)\leq \sup_{N\geq1}\frac{\E[z^{(N)}(\infty)]}{a}\leq \frac{1}{a}\Big(\lambda+\frac{\lambda}{1-\lambda}\Big),
}
This shows that to make $\sup_{N\geq1}\mathbb{P}(z^{(N)}(\infty)>a) < \epsilon$ for any $\epsilon > 0$, there exists appropriate choice
$a(\epsilon)$ not dependent on $N$. This shows that 
the sequence $\brac{z^{(N)}(\infty)}_{N}$ is tight.
The rest of the lemma now follows from the same line of arguments
as in the proof of the second statement of Theorem~\ref{thm:SA-JSQ_global_stability_limit_interchange}.
\end{proof}

The proof of Proposiion~\ref{thm:mf_rp} follows
directly from Lemma~\ref{lem:RP_global_stability}
by noting that $z^{(N)}(\infty)=Z^{(N)}(\infty)/N$.\qed

\section{Proof of Proposition~\ref{prop:tightness}}
\label{app:tightness}
Fix any $\epsilon >0$ and $l \geq 1$.
Using Markov inequality, we obtain
\eq{
\mathbb{P} \Big(\max_{j\in[M]}\sum_{i \geq l}x_{i,j}^{(N)}(\infty)>\epsilon\Big)
\leq\frac{\E\sbrac{\max_{j\in[M]}\sum_{i \geq l}x_{i,j}^{(N)}(\infty)}}{\epsilon}\leq \frac{1}{\epsilon}\E\sbrac{\sum_{j\in[M]}\sum_{i \geq l}x_{i,j}^{(N)}(\infty)}.
}
Since $\brac{x_{i,j}^{(N)}(\infty)}_i$ is a sequence non-negative random variables
for each $j \in [M]$, using monotone convergence theorem we can interchange the sum and the expectation on the RHS.
Hence, we have
\eqn{
\mathbb{P} \Big(\max_{j\in[M]}\sum_{i \geq l}x_{i,j}^{(N)}(\infty)>\epsilon\Big)\leq\frac{1}{\epsilon}\sum_{j\in [M]}\sum_{i\geq l}\E\sbrac{x_{i,j}^{(N)}(\infty)}=\frac{1}{\epsilon}\sum_{j\in [M]}\sum_{i\geq l}\mathbb{P} \sbrac{Q_{k,j}^{(N)}(\infty)\geq i},
\label{eq:tail_prob}
}
where the last equality follows from~\eqref{eqn:Realtion_Q_X}. Now, from Theorem~\ref{thm:tail_bound} we know
that for any $\theta \in [0,-\log \lambda)$ we have

\eq{
\sum_{j\in [M]}\sum_{i\geq l}\mathbb{P} \sbrac{Q_{k,j}^{(N)}(\infty)\geq i}\leq \sum_{j \in [M]}\sum_{i \geq l} C_j(\lambda,\theta)e^{-i\theta} = C(\theta) e^{-l\theta},
}
where $C(\theta)=\frac{(1-\lambda)}{(1-\lambda e^\theta)(1-e^{-\theta})}\sum_{j \in[M]}\frac{1}{\mu_j \gamma_j}$.
Since the RHS of the above inequality is not dependent on $N$,
using~\eqref{eq:tail_prob} we have

\begin{equation}
    \limsup_{N\to \infty} \mathbb{P} \Big(\max_{j\in[M]}\sum_{i \geq l}x_{i,j}^{(N)}(\infty)>\epsilon\Big) \leq C(\theta) e^{-l\theta}
\end{equation}
for all $\theta \in [0,-\log \lambda)$ and all $l \geq 1$. 
Hence, the condition of tightness given
by \eqref{eqn:l1_tight_3} is verified by 
fixing some $\theta \in (0, -\log \lambda)$
and letting $l \to \infty$.\qed

\section{Martingale Representation and Its Convergence}
\label{app:martingle_rep}
We first write the evolution of $x^{(N)}_{i,j}(t)$ for all $t\geq0$ in terms of number of arrivals and departures from system till time $t$, for $i\geq1$ and $j\in[M]$  as
\EQN{
\label{eqn:SA-JSQ_Flow}
x^{(N)}_{i,j}(t)=x^{(N)}_{i,j}(0)+ \frac{1}{N\gamma_j}A_{i,j}\Big(N \lambda \int_0^t p^{(N)}_{i-1,j}(\mathbf{x}^{(N)}(s))ds \Big)
-\frac{1}{N\gamma_j}D_{i,j}\Big(N\gamma_j\mu_j\int_0^t (x^{(N)}_{i,j}(s)-x^{(N)}_{i+1,j}(s))ds  \Big),
}
where $A_{i,j}$ and $D_{i,j}$ are mutually independent unit-rate Poisson processes and $p^{(N)}_{i-1,j}(\mathbf{x}^{(N)}(s))= \mathbbm{1}_{\{\mathbf{V}^{(N)}(s)\in R_{i,j}\}}$. Define for all $i\geq1$ and for all $j\in[M]$ 
\EQ{
\begin{aligned}
 M^{(A,N)}_{i,j}(t)&=A_{i,j}\Big(N \lambda \int_0^t p^{(N)}_{i-1,j}(\mathbf{x}^{(N)}(s))ds\Big)-N \lambda \int_0^t p^{(N)}_{i-1,j}(\mathbf{x}^{(N)}(s))ds,\\   
 M^{(D,N)}_{i,j}(t)&=D_{i,j}\Big(N \gamma_j\mu_j\int_0^t (x^{(N)}_{i,j}(s)-x^{(N)}_{i+1,j}(s))ds\Big)-N\gamma_j\mu_j\int_0^t (x^{(N)}_{i,j}(s)-x^{(N)}_{i+1,j}(s))ds.
\end{aligned}
}
We next show that the processes $M^{(A,N)}_{i,j}=(M^{(A,N)}_{i,j}(t):t\geq0)$ and $M^{(D,N)}_{i,j}=(M^{(D,N)}_{i,j}(t):t\geq0)$ are martingales with respect to the filtration $\mathbf{F}^{(N)}=\cbrac{F^{(N)}_t:t\geq0}$ augmented with all null sets, where $F^{N}_t=\bigcup_{j\in [M]}G_{t,j}^{(N)}$ with 
\eq{
\begin{aligned}
G_{t,j}^{(N)}=\bigcup_{i\geq1}\sigma\Big(
x^{(N)}_{i,j}(0),&A_{i,j}\Big(N \lambda \int_0^s p^{(N)}_{i-1,j}(\mathbf{x}^{(N)}(u))du \Big),\\
&D_{i,j}\Big(N\gamma_j\mu_j\int_0^s (x^{(N)}_{i,j}(u)-x^{(N)}_{i+1,j}(u))du \Big),0\le s\le t\Big).
\end{aligned}
}
\begin{lem}
\label{lem:predic_quadractic}
The processes $M^{(A,N)}_{i,j}$ and $M^{(D,N)}_{i,j}$, are square integrable $\mathbf{F}^{(N)}$-martingales for all $i\ge1$ and for all $j\in[M]$. Moreover, the predictable quadratic variation processes are given by
\EQ{
\begin{aligned}
\langle M^{(A,N)}_{i,j}\rangle(t)&=N \lambda \int_0^t p^{(N)}_{i-1,j}(\mathbf{x}^{(N)}(s))ds,\ i\geq1, \ j\in[M], \ t\geq0,\\
\langle M^{(D,N)}_{i,j}\rangle(t)&=N\gamma_j\mu_j\int_0^t (x^{(N)}_{i,j}(s)-x^{(N)}_{i+1,j}(s))ds,\ i\geq1, \ j\in[M], \ t\geq0.
\end{aligned}
}
\end{lem}
\begin{proof}
Let $I_{i,j}(t)=N \mu_j \gamma_j \int_0^t (x^{(N)}_{i,j}(s)-x^{(N)}_{i+1,j}(s))ds$ and $L_{i,j}(t)=N \lambda \int_0^t p^{(N)}_{i-1,j}(\mathbf{x}^{(N)}(s))ds$. Next, we prove that
\eq{
\begin{aligned}
\E(I_{i,j}(t))<\infty,& \ \E\Big(D_{i,j}(I_{i,j}(t))\Big)<\infty, \ i\geq1, \ j\in[M],\\ \E(L_{i,j}(t))<\infty,& \ \E\Big(A_{i,j}(L_{i,j}(t))\Big)<\infty, \ i\geq1, \ j\in[M].
\end{aligned}
}
The proof of square-integarble martingales and its corresponding predictable quadratic variation processes then
follows immediately using Lemma 3.2 from~\cite{Pang2007}.
Note that using crude inequality we can write
\eq{
\begin{aligned}
\E(I_{i,j}(t))&\leq \mu_jt\Big(N\gamma_j\E (x_{i,j}^{(N)}(0))+\E\Big(A_{i,j}(L_{i,j}(t))\Big) \Big)\\
&\leq \mu_j\gamma_j Nt+ \mu_jN\lambda t^2<\infty, \ t\geq0, \ i\geq1, \ j\in[M],
\end{aligned}
}
and
\eq{
\begin{aligned}
\E \Big[D_{i,j}(I_{i,j}(t))\Big]&\leq \E \Big[D_{i,j}\Big(\mu_jt\Big( N\gamma_j x_{i,j}^{(N)}(0)+A_{i,j}(L_{i,j}(t)) \Big)\Big)     \Big]\\
&=\E \cbrac{\E \Big[D_{i,j}\Big(\mu_jt\Big( N\gamma_j x_{i,j}^{(N)}(0)+A_{i,j}(L_{i,j}(t)) \Big)\Big)|A_{i,j}(L_{i,j}(t))\Big]}\\
&\leq \mu_jt(N\gamma_j+ N\lambda t)<\infty, t \geq0, \ i\geq1, \ j\in[M].
\end{aligned}
}
Similarly, we can show $\E(L_{i,j}(t))<\infty$ and $\E\Big(A_{i,j}(L_{i,j}(t))\Big)<\infty$.
\end{proof}

Now we can write the martingale representation of equation~\eqref{eqn:SA-JSQ_Flow} for all $t\geq0$, $i\geq1$, and $j\in[M]$ as  
\begin{multline}
\label{eqn:SA-JSQ_Mart_Rep_2}
x^{(N)}_{i,j}(t)=x^{(N)}_{i,j}(0)+ \frac{\lambda}{\gamma_j} \int_0^t \mathbbm{1}_{\{\mathbf{V}^{(N)}(s)\in R_{i,j}\}}ds  - \mu_j\int_0^t (x^{(N)}_{i,j}(s)-x^{(N)}_{i+1,j}(s))ds\\
+\frac{1}{N\gamma_j}(M^{(A,N)}_{i,j}(t)-M^{(D,N)}_{i,j}(t)).
\end{multline}
In below lemma, we prove that the martingale part in~\eqref{eqn:SA-JSQ_Mart_Rep_2} converges to $0$ with respect to $\ell_1$.
\begin{lem}
\label{lem:SA-JSQ_Mart}
Following convergence holds as $N\to \infty$
\eq{
\cbrac{\max_{j\in[M]}\frac{1}{N\gamma_j}\sum_{i\geq1}(|M^{(A,N)}_{i,j}(t)|+|M^{(D,N)}_{i,j}(t)|)}_{t\geq0} \Rightarrow 0. 
}
\end{lem}
\begin{proof}
The proof is similar to the proof of Proposition~4.3 in~\cite{Mukherjee2018} using Doob's inequality.
\end{proof}

\section{Proof of Lemma~\ref{lem:back_rate_thm1}}
To prove~\eqref{eqn:back_rate}, we consider the following linear optimisation problem

\eqn{
\label{eqn:Linear_Programme}
\begin{aligned}
\max \quad & \sum_{j\in[M]}\mu_j X^{(N,\Pi)}_{1,j}(t), \\
\textrm{s.t.} \quad & \sum_{j\in[M]} X^{(N,\Pi)}_{1,j}(t)\leq i,\\
  &0\leq X^{(N,\Pi)}_{1,j}(t) \leq N\gamma_j \ \forall j\in[M].\\
\end{aligned}
}
The first constraint in~\eqref{eqn:Linear_Programme} is true as the total number of busy servers  $\sum_{j\in[M]} X^{(N,\Pi)}_{1,j}(t)$ at time $t$ in $\mc M_N$ is always less than equal to the total number of customers in system that is $i$. The above optimisation problem has a unique maximum which is obtained in following way. First, note that the objective function is given as the weighted sum of $ X^{(N,\Pi)}_{1,j}(t)$ with weight $\mu_j$ for $j\in[M]$. We know that the maximum weight is $\mu_1$. Therefore, the maximum value that $X^{(N,\Pi)}_{1,1}(t)$ takes is $(i\land N\gamma_1 )$. Moreover, the second maximum weight is $\mu_2$, hence the maximum value that $X^{(N,\Pi)}_{1,2}(t)$ takes is $((i-N\gamma_1)_+\land N\gamma_2 )$. Proceeding in this way the maximum value that $X^{(N,\Pi)}_{1,j}(t)$ takes is $\brac{\brac{i-\sum_{i=1}^{j-1}N\gamma_i}_+ \land N\gamma_j}$ for $j\in[M]$, which completes the proof.  \qed

\section{Characterisation of Compact Sets and Tightness Criteria}
\label{sec:compact_sets}

In the lemma below, we characterise compact sets in the space $S$.
\begin{lem}
\label{lem:compact_sets_in_S}
A set $\mathbf{B}\subseteq S$ is relative compact in $S$ if and only if 
\eqn{
\label{eqn:relative_compact_set}
\lim_{l\to \infty}\sup_{\mathbf{y}\in \mathbf{B}}\max_{j\in[M]}\sum_{i\geq l}y_{i,j}=0.
}
\end{lem}
\begin{proof}
Suppose any $\mathbf{B}\subseteq S$ satisfying~\eqref{eqn:relative_compact_set}. To show $\mathbf{B}$ is relatively compact in $S$, we need to show that any sequence $(\mf y^{(N)})_{N\geq1}$ in $\mathbf{B}$ has a Cauchy subsequence. Since $S$ is complete under $\ell_1$, therefore the sequence $(\mf y^{(N)})_{N\geq1}$ has a convergent subsequence whose limit lies in $\mathbf{\bar{B}}$ which will complete the proof.

Next we show that the sequence $(\mf y^{(N)})_{N\geq1}$ has a Cauchy subsequence. Fix any $\epsilon>0$ and choose $l\geq1$ such that 
\eqn{
\label{eqn:relative_compact_set_2}
\max_{j\in [M]} \sum_{i \geq l}|y_{i,j}^{(N)}| < \frac{\epsilon}{4}, \ \forall N\geq1.
}
Now consider the sequence of first coordinates $(y_{1,j}^{(N)})_{N\geq1}$ for each $j\in[M]$. The sequence $(y_{1,j}^{(N)})_{N\geq1}$ lies in $[0,1]$. Therefore, by Bolzano-Wiestrass theorem it has a convergent subsequence $(y_{1,j}^{(N_k)})_{k\geq1}$. Moreover, along the indices $(N_k)_{k\geq1}$, the sequence $(y_{2,j}^{(N_k)})_{k\geq1}$ has a further convergent subsequence. Proceeding this way, we get a sequence of indices $(N_m)_{m\geq1}$ along which all first $l-1$ coordinates converges. This implies that there exists a $N^{\epsilon}\in \Nats$ such that 
\eqn{
\label{eqn:relative_compact_set_3}
\max_{j\in[M]} \sum_{i<l}|y_{i,j}^{(N)}-y_{i,j}^{(R)}|<\frac{\epsilon}{2}, \ \forall R,N\geq N^{\epsilon}.
}
Using~\eqref{eqn:relative_compact_set_2}, and~\eqref{eqn:relative_compact_set_3}, for all $R,N\geq N^{\epsilon}$ and for $\mf y^{(N)},\mf y^{(R)}\in \mf B$ we have
\eq{
\begin{aligned}
\norm{\mathbf{y}^{(N)}-\mathbf{y}^{(R)}}_1&=\max_{j\in[M]}\sum_{i\geq1}|y_{i,j}^{(N)}-y_{i,j}^{(R)}|\\
&\leq \max_{j\in [M]} \sum_{i < l}|y_{i,j}^{(N)}-y_{i,j}^{(R)}|+\max_{j\in [M]} \sum_{i \geq l}|y_{i,j}^{(N)}|+\max_{j\in [M]} \sum_{i \geq l}|y_{i,j}^{(R)}|\\
&<\epsilon,
\end{aligned}
}
along the sequence of indices $(N_m)_{m\geq1}$. This shows the existence of a Cauchy subsequence. Moreover, the limit point lies in $S$ follows from the completeness of $\ell_1$ space and the fact that $S$ is a closed subset of $\ell_1$.

Now for the only if part, let $\mathbf{B}$ be a relatively compact set in $S$. Assume that there exists a $\epsilon>0$ such that 
\eqn{
\label{eqn:comp_temp}
\lim_{l\to \infty}\sup_{\mathbf{y}\in \mathbf{B}}\max_{j\in[M]}\sum_{i\geq l}y_{i,j}>\epsilon.
}
This implies that for each $k\geq1$, there exists a $\mathbf{y}^{(k)}\in \mf B$, such that $\max_{j\in[M]}\sum_{i\geq l}y_{i,j}^{(k)}\geq\frac{\epsilon}{2}$. Therefore, if $\mathbf{y}^*$ be the limit of the sequence $(\mathbf{y}^{(k)})_{k\geq1}$, then~\eqref{eqn:comp_temp} implies that $\max_{j\in[M]}\sum_{i\geq l}y_{i,j}^{*}\geq\frac{\epsilon}{2}$ for all $l\geq1$, this leads to the contradiction that $\mathbf{y}^*\in \ell_1$.

\end{proof}
Next, we prove the criteria for a sequence to be tight in $S$.
\begin{lem}
\label{lem:tightness_condition}
A sequence $(\mathbf{y}^{(N)})_{N\geq1}$ of random elements in $S$ is tight iff for all $\epsilon>0$ we have
\eqn{
\label{eqn:l1_tight_1}
\lim_{l\to \infty} \limsup_{N\to \infty} \mathbb{P} \Big(\max_{j\in[M]}\sum_{i \geq l}y_{i,j}^{(N)}>\epsilon\Big)=0.
}
\end{lem}
\begin{proof}
For if part, we construct a relatively compact set $\mathbf{B}^{\epsilon}$ for any $\epsilon>0$ such that 
\eq{
\mathbb{P}\Big(\mathbf{y}^{(N)}\not\in \mathbf{B}^{\epsilon}\Big)<\epsilon, \ \forall N\in\Nats.
}
As $(\mathbf{y}^{(N)})_{N\geq1}$ satisfies~\eqref{eqn:l1_tight_1}, therefore there exists a $l(\epsilon)\geq1$ for all $\epsilon>0$ such that
\eq{
\limsup_{N\to\infty}\mathbb{P}\Big(\max_{j\in[M]}\sum_{i\geq l(\epsilon)}y_{i,j}^{(N)}>\epsilon\Big)<\epsilon,
}
and there exists a $N^{\epsilon}\geq1$ such that
\eq{
\mathbb{P}\Big(\max_{j\in[M]}\sum_{i\geq l(\epsilon)}y_{i,j}^{(N)}>\epsilon\Big)<\epsilon, \ \forall N>N^{\epsilon}.
}
Moreover, since $\mathbf{y}^{(1)},\dots,\mathbf{y}^{(N^{\epsilon})}$ are random elements of $S$, there exists $k(\epsilon)=\max\cbrac{l_1(\epsilon),\dots,l_{N^{\epsilon}}(\epsilon)}$ such that 
\eq{
\mathbb{P}\Big(\max_{j\in[M]}\sum_{i\geq k(\epsilon)}y_{i,j}^{(N)}>\epsilon\Big)<\epsilon, \ \forall N\in\Nats.
}
This implies that there exists an increasing sequence $(l(n))_{n\geq1}$ such that
\eq{
\mathbb{P}\Big(\max_{j\in[M]}\sum_{i\geq l(n)}y_{i,j}^{(N)}>\frac{\epsilon}{2^n}\Big)< \frac{\epsilon}{2^n}, \ \forall N\in\Nats.
}
Define 
\eq{
\mathbf{B}^{\epsilon}=\cbrac{\mathbf{y}\in S:\max_{j\in[M]}\sum_{i\geq l(n)}y_{i,j}\leq \frac{\epsilon}{2^n}, \forall n\geq1}.
}
From Lemma~\ref{lem:compact_sets_in_S}, the set $\mathbf{B}^{\epsilon}$ is relatively compact in $S$. Therefore, we have
\eq{
\mathbb{P}\Big(\mathbf{y}^{(N)}\not\in \mathbf{B}^{\epsilon}\Big)=\mathbb{P}\Big(\bigcup_{n\geq1}\cbrac{\max_{j\in[M]}\sum_{i\geq l(n)}y_{i,j}^{(N)}>\frac{\epsilon}{2^n}}\Big)\leq\sum_{n\geq1}\mathbb{P}\Big(\max_{j\in[M]}\sum_{i\geq l(n)}y_{i,j}^{(N)}>\frac{\epsilon}{2^n}\Big)<\epsilon,
}
where the first inequality follows from union bound. For only if part, assume that there exists a $\epsilon>0$ such that 
\eqn{
\label{eqn:l1_tight_2}
\lim_{l\to \infty} \limsup_{N\to \infty} \mathbb{P} \Big(\max_{j\in[M]}\sum_{i \geq l}y_{i,j}^{(N)}>\epsilon\Big)>\epsilon.
}
Since $(\mathbf{y}^{(N)})_{N\geq1}$ is tight in $S$, therefore there exists a convergent subsequence $(\mathbf{y}^{(N_k)})_{k\geq1}$ with limit $\mathbf{y}^*$. From~\eqref{eqn:l1_tight_2}, we can write
\eq{
\begin{aligned}
\epsilon<\lim_{l\to \infty} \limsup_{k\to \infty}\mathbb{P} \Big(\max_{j\in[M]}\sum_{i \geq l}y_{i,j}^{(N_k)}>\epsilon\Big)&\leq \lim_{l\to \infty} \limsup_{k\to \infty}\mathbb{P} \Big(\max_{j\in[M]}\sum_{i \geq l}y_{i,j}^{(N_k)}\geq\epsilon\Big)\\
&\leq \lim_{l\to \infty} \mathbb{P} \Big(\max_{j\in[M]}\sum_{i \geq l}y_{i,j}^*\geq\epsilon\Big),
\end{aligned}
}
where the last inequality follows from Portmanteau's theorem for closed set.
This leads to the contradiction that $\mathbf{y}^* \in \ell_1$.


\end{proof}

\section{Proof of Lemma~\ref{lem:relative_compact_SA-JSQ}}
\label{app:relative_completness_proof}
To prove relative compactness of the sequence $((\mathbf{x}^{(N)},\beta^{(N)}))_{N\geq1}$, we start with proving that 
for all finite time $t$ the system occupancy state lies in some compact set.
\begin{lem}
\label{lem:relative_3}
Assume $\mathbf{x}^{(N)}(0) \Rightarrow \mathbf{x}(0)\in S$, as $N\to \infty$. Then for any $T\geq0$, there exists a $L(T,\mathbf{x}(0))>2$ such that under the SA-JSQ policy, the probability that an arriving job joins a server with at-least $L(T,\mathbf{x}(0))-1$ active jobs upto time $T$ tends to $0$ as $N\to \infty$. 
\end{lem}
\begin{proof}
Let $A^{(N)}(t)$ denote the total number of arrivals up to time $t$. Arrivals are happening with rate $N\lambda$. Then for any $\epsilon>0$ we have 
\eq{
\mathbb{P}\Big(A^{(N)}(t)\geq (\lambda t+\epsilon)N \Big)\to 0, \ as \ N\to \infty.
}
Let $m^{(N)}(t)=\min_{j\in[M]}m^{(N)}_j(t)$, where $m^{(N)}_j(t)$ is the minimum queue length in the $j^{\textrm{th}}$ pool at time $t$ for $N^{\textrm{th}}$ system. The probability at which an arrival joins a server with at-least $L(T,\mathbf{x}(0))-1$ jobs during the time interval $[0,T]$ is given by
\begin{align}
&\mathbb{P}\Big(\cbrac{m^{(N)}(t)\geq L(T,\mathbf{x}(0))-1, \  t\in[0,T]} \cap \cbrac{an \ arrival \ occur \ at \ t}\Big)\nonumber\\
&\leq \mathbb{P}\Big(m^{(N)}(t)\geq L(T,\mathbf{x}(0))-1, \ t\in[0,T]\Big)\nonumber\\
&= \mathbb{P}\Big(X^{(N)}(t)\geq N(L(T,\mathbf{x}(0))-1), \ t\in[0,T]\Big),\label{eqn:lemma10_1}
\end{align}
where $X^{(N)}(t)$the is the total number of jobs in system at time $t$. From conservation of flow, we can write $X^{(N)}(t)=X^{(N)}(0)+A^{(N)}(t)-D^{(N)}(t)$, where $D^{(N)}(t)$ is the total number of departures from system till time $t$. Therefore, from~\eqref{eqn:lemma10_1}, we can write

\eqn{
\label{eqn:lemma10_2}
\begin{aligned}
&\mathbb{P}\Big(X^{(N)}(t)\geq N(L(T,\mathbf{x}(0))-1), \ t\in[0,T]\Big)\\
&=\mathbb{P}\Big(X^{(N)}(0)+A^{(N)}(t)-D^{(N)}(t)\geq N(L(T,\mathbf{x}(0))-1), \ t\in[0,T]\Big)\\
&\leq \mathbb{P}\Big(X^{(N)}(0)+A^{(N)}(t)\geq N(L(T,\mathbf{x}(0))-1), \ t\in[0,T]\Big)\\
&=\mathbb{P}\Big( M^{(N)}(t) \geq (L(T,\mathbf{x}(0))-1)-x^{(N)}(0)-\lambda t, \ t\in[0,T]\Big),
\end{aligned}
}

where $M^{(N)}(t)=\frac{A^{(N)}(t)-N\lambda t}{N}$, and $x^{(N)}(0)=\frac{X^{(N)}(0)}{N}$. Now observer that for all $\epsilon>0$ we have
\EQ{
\mathbb{P}\Big(\sup_{t\in[0,T]}|M^{(N)}(t)|>\epsilon\Big) \to 0, \ as \ N \to \infty.
}
Therefore, we have
\EQN{
\label{eqn:lemma10_3}
\mathbb{P}\Big(|M^{(N)}(t)|>1, \ t\in[0,T]\Big) \to 0, \ as \ N \to \infty.
}
From~\eqref{eqn:lemma10_1},~\eqref{eqn:lemma10_2},~\eqref{eqn:lemma10_3} and choosing $L(T,\mathbf{x}(0))>2+x(0)+\lambda T$, where $x(0)=\sum_{j\in[M]}\sum_{i\geq 1}\gamma_j x_{i,j}(0)$, we have
\eq{
\mathbb{P}\Big(\cbrac{m^{(N)}(t)\geq L(T,\mathbf{x}(0))-1, \  t\in[0,T]} \cap \cbrac{an \ arrival \ occur \ at \ t}\Big) \to 0, \ as \ N\to\infty.
}

\end{proof}
We now prove the relative compactness of the sequence $((\mathbf{x}^{(N)},\beta^{(N)}))_{N\geq1}$ by showing the relative compactness of individual component. Note that the space $\mathbf{E}$ is compact. Therefore, relative compactness of the sequence $(\beta^{(N)})_{N\geq1}$ follows from Prohorov's theorem~\cite{Billingsley2013}. To prove relative compactness of the sequence $(\mathbf{x}^{(N)})_{N\geq1}$, we need to verify following conditions.
\begin{enumerate}
    \item For every $\eta>0$ and rational $t\geq0$, there exists a compact set $\mathbf{B}_{\eta,t}\subset S$ such that
    \eqn{
    \label{eqn:cond_rel_1}
    \liminf_{N\to \infty} \mathbb{P}(\mathbf{x}^{(N)}(t)\in \mathbf{B}_{\eta,t})\geq1-\eta}
    \item For every $\eta>0$ and for $T>0$ there exists a $\delta>0$ and a finite partition $\cbrac{t_1,t_2,\dots,t_n}$ of $[0,T]$ with $\min_{l\in[n]} |t_l-t_{l-1}|>\delta$ such that 
    \eqn{
    \label{eqn:cond_rel_2}
    \limsup_{N\to \infty}\mathbb{P}\Big(\max_{l\in[n]} \sup_{s,t \in [t_{l-1},t_l)}\norm{\mathbf{x}^{(N)}(s)-\mathbf{x}^{(N)}(t)}_1 \geq \eta \Big) <\eta.}
\end{enumerate}
We first prove condition~\eqref{eqn:cond_rel_1}.
From Lemma~\ref{lem:relative_3}, for any fix $t\geq0$ we have 
\eq{\lim_{N\to \infty} \mathbb{P}\Big( x_{i,j}^{(N)}(t) \leq x_{i,j}^{(N)}(0), \ \forall i\geq L(t,\mathbf{x}(0)),j\in[M]\Big)=1.}
 Now it can be easily verified that the sequence $(\mathbf{x}^{(N)}(0))_{N\geq1}$ is tight in $S$. Therefore, using previous condition we can write for any $\epsilon>0$
\eqn{
\label{eqn:relative_1}
\lim_{l\to \infty} \limsup_{N \to \infty}\mathbb{P}\Big(\max_{j\in[M]}\sum_{i\geq l}x_{i,j}^{(N)}(t)>\epsilon\Big)\leq \lim_{l\to \infty} \limsup_{N \to \infty}\mathbb{P}\Big(\max_{j\in[M]}\sum_{i\geq l}x_{i,j}^{(N)}(0)>\epsilon\Big)=0.
}
Hence, from Lemma~\ref{lem:tightness_condition} the sequence $(\mathbf{x}^{(N)}(t))_{N\geq1}$ is tight in $S$. This implies that the condition~\eqref{eqn:cond_rel_1} is satisfied. 
Next for any $t_1<t_2$, we consider
\eq{
\begin{aligned}
&|x_{i,j}^{(N)}(t_1)-x_{i,j}^{(N)}(t_2)|\leq \frac{\lambda}{\gamma_j}\beta^{(N)}([t_1,t_2]\times R_{i,j})+\mu_j\int_{t_1}^{t_2} (x^{(N)}_{i,j}(s)-x^{(N)}_{i+1,j}(s))ds\\
&+\frac{1}{N\gamma_j}|M^{(A,N)}_{i,j}(t_1)-M^{(D,N)}_{i,j}(t_1) - M^{(A,N)}_{i,j}(t_2)+M^{(D,N)}_{i,j}(t_2)|+o(1), \ i\geq1,j\in[M].
\end{aligned}
}
Using the above equation we can write $\ell_1$ distance between $\mathbf{x}^{(N)}(t_1)$ and $\mathbf{x}^{(N)}(t_2)$ as
\eqn{
\label{eqn:relative_2}
\begin{aligned}
&\norm{\mathbf{x}^{(N)}(t_1)-\mathbf{x}^{(N)}(t_2)}_{1}\\
&\leq\max_{j\in[M]}\sum_{i\geq1}\frac{\lambda}{\gamma_j}\beta^{(N)}([t_1,t_2]\times R_{i,j})+\max_{j\in[M]}\sum_{i\geq1}\int_{t_1}^{t_2} \mu_j(x^{(N)}_{i,j}(s)-x^{(N)}_{i+1,j}(s))ds\\
&+\max_{j\in[M]}\frac{1}{N\gamma_j}\sum_{i\geq1}|M^{(A,N)}_{i,j}(t_1)-M^{(D,N)}_{i,j}(t_1) - M^{(A,N)}_{i,j}(t_2)+M^{(D,N)}_{i,j}(t_2)|+o(1)\\
&\leq \frac{\lambda}{\gamma_{min}}(t_2-t_1)+\max_{j\in[M]}\int^{t_2}_{t_1}\mu_jx_{1,j}^{(N)}(s)ds+\max_{j\in[M]}\frac{1}{N\gamma_j}\sum_{i\geq1}|M^{(A,N)}_{i,j}(t_1)-M^{(D,N)}_{i,j}(t_1)(s)ds\\
&- M^{(A,N)}_{i,j}(t_2)+M^{(D,N)}_{i,j}(t_2)|+o(1)\\
&\leq (\frac{\lambda}{\gamma_{min}}+\mu_1)(t_2-t_1)+\max_{j\in[M]}\frac{1}{N\gamma_j}\sum_{i\geq1}|M^{(A,N)}_{i,j}(t_1)-M^{(D,N)}_{i,j}(t_1)- M^{(A,N)}_{i,j}(t_2)+M^{(D,N)}_{i,j}(t_2)|+o(1),
\end{aligned}
}
where $\gamma_{min}=\min_{j\in[M]}\gamma_j$.
From Lemma~\ref{lem:SA-JSQ_Mart}, the martingale part in~\eqref{eqn:relative_2} converges to $0$ as $N\to \infty$.
Moreover, from~\eqref{eqn:relative_2}, it implies that for any finite partition $\cbrac{t_1,t_2,\dots,t_n}$ of $[0,T]$ with $\min_{l\in[n]}|t_l-t_{l-1}|>\delta$, we have 
\eq{
\max_{l\in[n]} \sup_{s,t \in [t_{l-1},t_l)}\norm{\mathbf{x}^{(N)}(s)-\mathbf{x}^{(N)}(t)}_1\leq (\frac{\lambda}{\gamma_{min}}+\mu_1)\max_{l\in[n]}(t_l-t_{l-1})+\zeta^{(N)},
}
where $\mathbbm{P}(\zeta^{(N)}>\frac{\eta}{2})<\eta$ for all sufficiently large $N$. Take $\delta=\eta/(4(\frac{\lambda}{\gamma_{min}}+\mu_1))$ and any partition with $\max_{l\in[n]}(t_l-t_{l-1})<\eta/(2(\frac{\lambda}{\gamma_{min}}+\mu_1))$ and $\min_{l\in[n]}|t_l-t_{l-1}|>\delta$, we have
\eq{
\max_{l\in[n]} \sup_{s,t \in [t_{l-1},t_l)}\norm{\mathbf{x}^{(N)}(s)-\mathbf{x}^{(N)}(t)}_1\leq\eta,
}
on the event $\{\zeta^{(N)}\leq\frac{\eta}{2}\}$. Therefore, for sufficiently large $N$ we obtain
\eq{
\mathbb{P}\Big( \max_{l\in[n]} \sup_{s,t \in [t_{l-1},t_l)}\norm{\mathbf{x}^{(N)}(s)-\mathbf{x}^{(N)}(t)}_1\geq\eta\Big)\leq\mathbb{P}(\zeta^{(N)}>\frac{\eta}{2})<\eta.
}
Hence the condition~\eqref{eqn:cond_rel_2} is satisfied.
Next, we show that the limit $(\mathbf{x},\beta)$ of any convergent subsequence of the sequence $((\mathbf{x}^{(N)},\beta^{(N)}))_{N\geq1}$ satisfies~\eqref{eqn:SA-JSQ_random_rep_limit}. We first show that the right side of~\eqref{eqn:outline_3} is a continuous map and then the result follows from an application of continuous mapping theorem. Consider,
\eq{
\begin{aligned}
W_{i,j}\Big(\mathbf{x}^{(N)}(t),\beta^{(N)},\mathbf{x}^{(N)}(0),\mathbf{m}^{(N)}\Big)(t)&= x^{(N)}_{i,j}(0)+ \frac{\lambda}{\gamma_j}\beta^{(N)}([0,t]\times R_{i,j})\\   &- \mu_j\int_0^t (x^{(N)}_{i,j}(s)-x^{(N)}_{i+1,j}(s))ds+ m_{i,j}^{(N)}(t) \ i\geq1, \ j\in[M],
\end{aligned}
}
where $m_{i,j}^{(N)}(t)=\frac{1}{N\gamma_j}(M^{(A,N)}_{i,j}(t)-M^{(D,N)}_{i,j}(t))$. Now we next prove that the map $\mathbf{W}=(W_{i,j})_{i,j}$ is continuous.
First, assume that the sequence $((\mathbf{x}^{(N)},\mathbf{m}^{N}))_{N\geq1}$ converges to $(\mathbf{x},\mathbf{m})$, then there exists a $N_1\in\Nats$ such that $\sup_{t\in[0,T]}\norm{\mathbf{x}^{(N)}(t)-\mathbf{x}(t)}_1<\epsilon/(4T)$ for all $N\geq N_1$. Therefore, we can write
\eq{
\sup_{t\in[0,T]}\int^t_0 |x_{1,j}^{(N)}(t)-x_{1,j}(t)|ds\leq T\sup_{t\in[0,T]}\norm{\mathbf{x}^{(N)}(t)-\mathbf{x}(t)}_1<\epsilon/4.
}
Also, there exists a $N_2\in\Nats$ such that $\sup_{t\in[0,T]}\norm{\mathbf{m}^{(N)}(t)-\mathbf{m}(t)}_1<\epsilon/4$.
Second, assume that the sequence $(\mathbf{x}^{(N)}(0))$ converges to $\mathbf{x}(0)$ with respect to $\ell_1$. Therefore, there exists a $N_3\in \Nats$ such that $\norm{\mathbf{x}^{(N)}(0)-\mathbf{x}(0)}_1<\epsilon/4$. Now we claim that there exists a $N_4 \in \Nats$ such that
\eqn{
\label{eqn:relative_3}
\frac{\lambda}{\gamma_j} \max_{j\in[M]}\sum_{i\geq1}|\beta^{(N)}([0,T]\times R_{i,j})-\beta([0,T]\times R_{i,j})|\leq \epsilon/4.
}
The equation~\eqref{eqn:relative_3} implies that the convergence of the sequence $((\beta^{(N)}([0,t]\times R_{i,j}))_{i,j})_{N\geq1}$ for any $t\geq0$ is $\ell_1$ convergence. However, we know only the weak convergence of the sequence of measures $(\beta^{(N)})_{N\geq1}$, which does not directly implies~\eqref{eqn:relative_3}. Therefore, we show using weak convergence of the sequence $(\beta^{(N)})_{N\geq1}$ that~\eqref{eqn:relative_3} is indeed true in our case. Since $\mf x(0)\in S$, there exists a $m_j'(\mf x(0))$ for all $j\in[M]$ such that $x_{i,j}(0)<1$ for all $i\geq m_j'(\mf x(0))$. Furthermore, from Lemma~\ref{lem:relative_3} we can write
\eq{\lim_{N\to \infty} \mathbb{P}\Big( \sup_{t\in[0,T]}x_{i,j}^{(N)}(t) \leq x_{i,j}^{(N)}(0), \ \forall i\geq L(T,\mathbf{x}(0)),j\in[M]\Big)=1.}
Therefore, for $N'=\max\cbrac{m_j'(\mf x(0)),L(T,\mathbf{x}(0))}$ we have
\eqn{
\label{eqn:relative_4}
\lim_{N\to \infty} \mathbb{P}\Big( \sup_{t\in[0,T]}x_{i,j}^{(N)}(t) \leq1, \ \forall i\geq N',j\in[M]\Big)=1.}
Using~\eqref{eqn:relative_4} and~\eqref{eqn:R_1} we get
\eqn{
\label{eqn:relative_5}
\lim_{N\to \infty}\max_{j\in[M]}\sum_{i\geq N'}\beta^{(N)}([0,T]\times R_{i,j}) =\max_{j\in[M]}\sum_{i\geq N'}\beta([0,T]\times R_{i,j})=0.
}
Also, weak convergence of $(\beta^{(N)})_N$ implies that
\eqn{
\label{eqn:relative_6}
\lim_{N\to \infty}\max_{j\in[M]}\sum_{i< N'}\beta^{(N)}([0,T]\times R_{i,j}) =\max_{j\in[M]}\sum_{i< N'}\beta([0,T]\times R_{i,j}).
}
Hence, using~\eqref{eqn:relative_5} and~\eqref{eqn:relative_6} we get the desired result that is~\eqref{eqn:relative_3}. Let $\bar{N}=\max\cbrac{N_1,N_2,N_3,N_4}$, then we have 
\eq{
\sup_{t\in[0,T]}\norm{\mathbf{W}\Big(\mathbf{x}^{(N)}(t),\beta^{(N)},\mathbf{x}^{(N)}(0),\mathbf{m}^{(N)} \Big)-\mathbf{W}\Big(\mathbf{x}(t),\beta,\mathbf{x}(0),\mathbf{m} \Big)}_1(t)<\epsilon.
}
This shows that the map $\mathbf{W}$ is continuous, which completes the proof.


\section{Proof of Lemma~\ref{lem:M-dim-Markov}}
To prove positive recurrent of the Markov chain $\mf U$, we consider the Lyapunov function $f:\mb{Z}_+^K \to [0,\infty)$ as 
\eq{
f( \mf R )=\sum_{j\in [K]}R_j, \ \mf R \in \mb{Z}_+^K.
}
Moreover, from~\eqref{eqn:k-dim-chain} we can write the generator $G_{\mf U}$ acting on the function $f$ for a state $\mf R \in\mb{Z}_+^K $ as
\eq{
G_{\mf U} f(\mf R)=\sum_{j\in[K]}\nu_j - \lambda\sum_{j\in[K]} \indic{0=R_{1}=\dots=R_{j-1}< R_{j}}.
}
Define the set $B=\cbrac{\mf S \in \mb{Z}_+^K:S_i=0 \forall i\in[K]}$. Observe that if $\sum_{j\in[K]}\nu_j<\lambda$, then $G_{\mf U} f(\mf R)<0$ for all $\mf R \in B^c$, where $B^c$ is the compliment of the set $B$ and is given by
\eq{
B^c=\cbrac{\mf S \in \mb{Z}_+^K: \mf S \in \bigcup_{i=1}^K \cbrac{0=S_{1}=\dots=S_{i-1}< S_{i}} }.
}
Otherwise $G_{\mf U} f(\mf R)<\sum_{j\in[K]}\nu_j$ for $\mf R\in B$.
Hence, using the Foster-Lyapunov criterion for positive recurrence from~\cite{Kelly_book},
we conclude that the chain $\mf U$
is positive recurrent if $\sum_{j\in[K]}\nu_j<\lambda$. 

Next, we proceed to prove that the unique stationary distribution $\pi$ of the chain $\mf U$ satisfies~\eqref{eqn:stat_prop_finite_M}.
We prove this using induction on $K$.
First observer that $\pi(U_1>0)=\nu_1/\lambda$ if the birth-death process corresponding to the component $U_1$ is stable. This proves the base case, i.e., $K=1$.
Now for the induction hypothesis assume that~\eqref{eqn:stat_prop_finite_M} is true for all $j\in[K-1]$. Note that the positive recurrence of the chain $\mf U$ implies that $\E f(\mf U) < \infty$, where the expectation is with respect to the stationary distribution $\pi$. Therefore, we can set the steady-state expected drift to $0$ which gives
\eqn{
\label{eqn:steady_state_drift}
\E[G_{\mf U} f(\mf U)]=\sum_{j\in[K]}\nu_j -\lambda \sum_{j\in[K]} \pi(0=U_1=\dots=U_{j-1}< U_{j})=0.
}
Hence, using the induction hypothesis for all $j\in[K-1]$ in the above expression we get
\eq{
\pi(0=U_1=\dots=U_{K-1}< U_{K})=\frac{\nu_K}{\lambda},
} 
which completes the proof for~\eqref{eqn:stat_prop_finite_M}.

Now we show that the condition $\sum_{j\in[K]}\nu_j<\lambda$ is necessary for the chain $\mf U$ to be positive recurrent. Suppose the chain $\mf U$ is positive recurrent and, therefore, it has a unique stationary measure $\pi$. From~\eqref{eqn:steady_state_drift} it is clear that $\pi$ must satisfy~\eqref{eqn:stat_prop_finite_M}. Therefore, we can write
\eq{
\pi\cbrac{U_i=0 \ \forall i\in[K] }= 1- \sum_{i=1}^K \pi\cbrac{0=U_1=\dots=U_{i-1}< U_{i}}=1-\sum_{i=1}^K \frac{\nu_i}{\lambda}>0,
}
where the last inequality follows since the chain $\mf U$ is positive recurrent. 
\qed

\section{Proof of Theorem~\ref{thm:SA-JSQ_global_stability_limit_interchange}   }
\label{app:global_stability}
{\bf Proof of Theorem~\ref{thm:SA-JSQ_global_stability_limit_interchange}.(\ref{global_stab_sa_jsq})}:
To prove this, we 
need the following lemma
which extends 
the monotonicity of the process
$\mf x^{(N)}$ for finite $N$
to the monotonicity of the limiting
process $\mf x$.

\begin{lem}
\label{lem:SA-JSQ_ODE_Solution_Monotone}
Let ${\mf x}(\cdot,\mf u)=({\mf x}(t,\mf u), t\geq 0)$
denote a solution to~\eqref{eqn:SA-JSQ_fluid_process} with $\mathbf{x}(0)=\mathbf{u}\in S$. 
Then, for any $\mf u, \mf v \in S$ satisfying $\mf u \leq \mf v$
we have ${\mf x}(t, \mf u) \leq {\mf x}(t, \mf v)$ for all $t \geq 0$.
\end{lem}
\begin{proof}
First, note that for any $\mathbf{u} \in S$, there exists a sequence $\brac{\mathbf{u}^{(N)}}_{N\geq 1}$ with $\mathbf{u}^{(N)}=(u^{(N)}_{i,j}, i\geq 1, j \in [M])\in S\cap S^{(N)}$ such that $\norm{\mathbf{u}^{(N)} - \mathbf{u}}_1 \to 0$ as $N \to \infty$. 
We can simply construct such a sequence by 
setting $u^{(N)}_{i,j}=\frac{\floor{u_{i,j}N \gamma_j}}{N\gamma_j}$ for each $i \geq 1$ and $j \in [M]$. This construction
also satisfies the property that if $\mf u, \mf v \in S$
are such that $\mf u \leq \mf v$ and
if the sequences $(\mf u^{(N)})_N$
and $(\mf v^{(N)})_N$ are constructed from their corresponding limits $\mf u$ and $\mf v$
as described above, 
then $\mf u^{(N)} \leq 
\mf v^{(N)}$ for all $N$. 

Let $\mf x^{(N)}(\cdot,\mf u^{(N)})=(\mf x^{(N)}(t,\mf u^{(N)}), t\geq 0)$ denote the process $\mf x^{(N)}$
started at $\mf x^{(N)}(0)=\mf u^{(N)}$.
Then by Theorem~\ref{thm:SA-JSQ_monotone}
we have that 

\begin{equation}
    \mf x^{(N)}(t, \mf u^{(N)}) \leq \mf x^{(N)}(t, \mf v^{(N)}), \forall t \geq 0.
\end{equation}
Now letting $N \to \infty$ and applying
Theorem~\ref{thm:SA-JSQ_Process_level_convergence}
gives the desired result.
\end{proof}

For $\mf u \in S$, we define $v_{n,j}(t,\mf u)=\sum_{i\geq n}x_{i,j}(t,\mf u)$ and $v_{n,j}(\mf u)=\sum_{i\geq n}u_{i,j}$
for each $n\geq1$ and $j\in[M]$.
Furthermore, let $v_n(t,\mf u)=\sum_{j\in[M]} \gamma_j v_{n,j}(t,\mf u)$ and
$v_n(\mf u)=\sum_{j\in[M]} \gamma_j v_{n,j}(\mf u)$ for each $n\geq1$ and $\mf u \in S$.

 \begin{lem}
 For $\mf u\in S$ let $ \mf x(\mf u,\cdot)$
 denote a solution of~\eqref{eqn:SA-JSQ_fluid_process}
 in $S$.
 Then for all $t\geq 0$ we have 
 \EQN{
 \label{eqn:v_n}
 \frac{dv_n(t,\mf u)}{dt}=\lambda \sum_{j\in[M]}\sum_{i\geq n}p_{i-1,j}(\mf x(t,\mf u))-\sum_{j\in[M]}\mu_j \gamma_j x_{n,j}(t,\mf u), \ \forall n\geq1.
 }
 In particular, we have 
 \EQN{\label{eqn:v_1}
 \frac{dv_1(t,\mf u)}{dt}=\lambda -\sum_{j\in[M]}\mu_j \gamma_j x_{1,j}(t,\mf u).
 }
\end{lem}
\begin{proof}
We can write the differential form of~\eqref{eqn:SA-JSQ_fluid_process} as
\eqn{
\frac{dx_{i,j}(t,\mf u)}{dt}=\frac{\lambda}{\gamma_j}p_{i-1,j}(\mf x(t,\mf u))-\mu_j (x_{i,j}(t,\mf u)-x_{i+1,j}(t,\mf u)).
}
Multiplying the above by $\gamma_j$ and summing first over $i\geq n$ and then over $j\in [M]$ we obtain~\eqref{eqn:v_n}.
Equation~\eqref{eqn:v_1} 
follows from~\eqref{eqn:v_n} by using 
P\ref{p2}.
\end{proof}

From Lemma~\ref{lem:SA-JSQ_ODE_Solution_Monotone} it follows
that for any $\mathbf{x}(0)\in S$ and any $t \geq 0$
we have
\EQN{
\label{eqn:sandwitch_sol}
\mathbf{x}(t,\min(\mathbf{x}(0),\mathbf{x}^*))\leq\mathbf{x}(t,\mathbf{x}(0))\leq \mathbf{x}(t,\max(\mathbf{x}(0),\mathbf{x}^*)),
}
where $\min(\mf u,\mf v)$ with $\mf u, \mf v \in S$
is defined by taking the component-wise minimum.
Hence, to prove global stability
it is sufficient to show that 
the (component-wise) convergence $\mf x(t,\mf x(0)) \to \mf x^*$ holds for
initial states satisfying either of the following two conditions:
(i) $\mf x(0) \geq \mf x^*$ and (ii) $\mf x(0) \leq \mf x^*$.

To prove the above, we first show that 
for any solution $\mf x(\cdot,\mf x(0)) \in S$, $v_n(t,\mf x(0))$ is uniformly bounded in $t$ for all $n \geq 1$.
Consider the case when $\mf x(0) \geq \mf x^*$. 
From Lemma~\ref{lem:SA-JSQ_ODE_Solution_Monotone}
it follows that for $\mf x(0) \geq \mf x^*$, we have
$\mf x(t,\mf x(0)) \geq \mf x^*$ for all $t \geq 0$. Therefore, we can write
\eq{
\sum_{j\in[M]}\gamma_j \mu_j x_{1,j}(t,\mf x(0))\geq \sum_{j\in[M]}\gamma_j \mu_j x_{1,j}^*=\lambda,
}
where the last equality follows from~\eqref{eqn:dep_arri_rate}.
Hence, from~\eqref{eqn:v_1} we have $\frac{dv_1(t,\mf x(0))}{dt}\leq 0$
from which it follows that
$0\leq v_1(t,\mf x(0))\leq v_1(\mf x(0))$ for all $t\geq0$. 
Since the sequence $\brac{v_n(t,\mf x(0))}_{n\geq 1}$ is non-increasing, we have $0\leq v_n(t,\mf x(0))\leq v_1(\mf x(0))$ for all $n\geq1$ and for all $t\geq0$. This proves that $v_n(t,\mf x(0))$ is uniformly bounded in $t$ for each $n\geq1$ if $\mf x(0) \geq \mf x^*$.
Now consider the case $\mf x(0)\leq \mf x^*$. From Lemma~\ref{lem:SA-JSQ_ODE_Solution_Monotone}
it follows that for $\mf x(0) \leq \mf x^*$, we have
$\mf x(t,\mf x(0)) \leq \mf x^*$ for all $t \geq 0$. Therefore, we have $v_1(t,\mf x(0))\leq v_1(\mf x^*)$ for all $t\geq0$. This shows that the component $v_n(t,\mf x(0))$ is uniformly bounded in $t$ for each $n\geq1$ for $\mf x(0)\leq \mf x^*$.

Since  $v_n(t,\mf x(0))$ is uniformly bounded in $t$, the convergence $x_{i,j}(t,\mf x(0)) \to x_{i,j}^*$  for all $i\geq 1$ and for all $j\in [M]$ will follow from 
\eqn{
\label{eqn:int_finite}
\int_0^{\infty} (x_{i,j}(t,\mf x(0)) - x_{i,j}^*) dt < \infty, \ \forall j\in[M], \ \forall i\geq1,
}
for the case $\mf x(0)\geq \mf x^*$ and from
\eqn{
\label{eqn:int_finite_2}
\int_0^{\infty} (x_{i,j}^* - x_{i,j}(t,\mf x(0))) dt < \infty, \ j\in[M], \ i\geq1,
}
for the case $\mf x(0)\leq \mf x^*$. We now prove~\eqref{eqn:int_finite} to show convergence for the case $\mf x(0)\geq \mf x^*$; the proof of other case follows similarly.
To show~\eqref{eqn:int_finite} it is sufficient to prove that
\eqn{
\label{eqn:int_finite_2}
\int_0^{\infty} 
\sum_{j\in [M]} \mu_j \gamma_j(x_{i,j}(t,\mf x(0)) - x_{i,j}^*) dt < \infty, \ \forall i\geq1.
}

For $i=1$, we can write~\eqref{eqn:int_finite_2} as
\eq{
\begin{aligned}
\int_0^{\tau} 
(\sum_{j\in [M]} \mu_j \gamma_jx_{1,j}(t,\mf x(0)) - \sum_{j\in [M]} \mu_j \gamma_jx_{1,j}^*) dt&=\int_0^{\tau} 
(\sum_{j\in [M]} \mu_j \gamma_jx_{1,j}(t,\mf x(0)) -\lambda)\\
&=-  \int_0^{\tau} \frac{dv_1(t,\mf x(0))}{dt}dt\\
&=v_1(\mf x(0))-v_1(\tau,\mf x(0))\\
&\leq v_1(\mf x(0)),
\end{aligned}
}
where the first equality follows from~\eqref{eqn:dep_arri_rate}, second equality follows from~\eqref{eqn:v_1}, and the last inequality follows as $v_1(t,\mf x(0))$ is uniformly bounded in $t$.
Since the right hand side is bounded by a constant for all $\tau$ , the integral on the left hand side must converge as $\tau \to \infty$. This shows that $x_{1,j}(t,\mf x(0)) \to x_{1,j}^*$ for all $j\in[M]$ as $t\to \infty$.

Now for $i=2$, we can write~\eqref{eqn:int_finite_2} as 
\begin{align}
\int_0^{\tau} (\sum_{j\in [M]} \mu_j \gamma_jx_{2,j}(t,\mf x(0)) -\sum_{j\in [M]} \mu_j \gamma_jx_{2,j}^*)dt&=\int_0^{\tau} \sum_{j\in [M]} \mu_j \gamma_jx_{2,j}(t,\mf x(0)) dt \nonumber\\
&=-  \int_0^{\tau} \Big(\frac{dv_2(t,\mf x(0))}{dt}- \lambda \sum_{j\in[M]}\sum_{i\geq2}p_{i-1,j} (\mf x(t,\mf x(0)))\Big)dt, \nonumber\\
&=\lambda \int_0^{\tau}\sum_{j\in[M]}\sum_{i\geq2}p_{i-1,j} (\mf x(t,\mf x(0)))dt + v_1(\mf x(0)) - v_1(\tau,\mf x(0)) \nonumber\\
&\leq \lambda \int_0^{\tau}\sum_{j\in[M]}\sum_{i\geq2}p_{i-1,j} (\mf x(t,\mf x(0)))dt + v_1(\mf x(0)) \label{eqn:second_comp},
\end{align}
where the first equality follows as $x_{2,j}^*=0$ for all $j\in [M]$ and the second equality follows from~\eqref{eqn:v_n}
for $n=2$. 
From the convergence $x_{1,j}(t,\mf x(0))\to x^*_{1,j}$ for all $j\in [M]$, we know that for all $\epsilon>0$ there exists a $t^{\epsilon}>0$ such that $(x_{1,j}(t,\mf x(0))-x_{1,j}^*)< \epsilon$ for all $t\geq t^{\epsilon} $ and for all $j \in [M]$. 
Furthermore, from~\eqref{eqn:SA-JSQ_Fixed_Point} it follows that for $\lambda <1$ there exists $j \in [M]$
such that $x^*_{1,j} <1$. Hence, we can choose
$\epsilon < 1-x^*_{1,j}$ which yields
%
\eq{
x_{1,j}(t,\mf x(0))<\epsilon +x_{1,j}^*<1, \ \forall t\geq t_0^{\epsilon}.}
%
Hence, by properties~P\ref{p2}, P\ref{p5}, and
P\ref{p6} of the fluid limit we have $\sum_{j\in[M]}\sum_{i\geq2}p_{i-1,j} (\mf x(t,\mf x(0)))=0$ for all $t\geq t_0^{\epsilon}$.
Therefore, for all $\tau \geq t^{\epsilon}$, we can write~\eqref{eqn:second_comp} as 
\begin{align}
\int_0^{\tau} \sum_{j\in [M]} \mu_j \gamma_jx_{2,j}(t,\mf x(0)) dt 
& \leq \lambda \int_0^{t^{\epsilon} }\sum_{j\in[M]}\sum_{i\geq2}p_{i-1,j} (\mf x(t,\mf x(0)))dt +v_1(\mf x(0))\leq \lambda t^{\epsilon}+v_1(\mf x(0))\nonumber,
\end{align}
where in the last inequality we have used the fact that $\sum_{j\in[M]}\sum_{i\geq2}p_{i-1,j} (\mf x(t,\mf x(0)))\leq1$ from P\ref{p2}.  
Since the right hand side of the above expression is independent of $\tau$, the left hand side must converge as $\tau \to \infty$. This shows that $x_{2,j}(t,\mf x(0)) \to x_{2,j}^*=0$ as $t\to \infty$ for all $j\in[M]$. 

Finally, since $\mf x(t,\mf x(0)) \in S$ we have $x_{2,j}(t,\mf x(0))\geq x_{i,j}(t,\mf x(0))\geq 0$ for all $i\geq 3$, $j\in[M]$ and for all $t\geq 0$. Hence, $x_{i,j}(t,\mf x(0)) \to 0$ as $t \to \infty$ for all $i\geq 3$ and $j\in[M]$.




{\bf Proof of Theorem~\ref{thm:SA-JSQ_global_stability_limit_interchange}.(\ref{limit_interchange})}:
We first
recall from Proposition~\ref{prop:tightness}
that the sequence $(\mf x^{(N)}(\infty))_N$
is tight in $S$ under the $\ell_1$-norm.
Hence, by Prohorov's theorem,
the sequence has convergent subsequences
with limits in $S$. Thus, it suffices to
show that all convergent
subsequences has the same limit point $\mf x^*$.
Let $(\mf x^{(N_k)}(\infty))_{k}$ be any such convergent subsequence of the sequence $(\mf x^{(N)}(\infty))_N$ with limit point $\mf x^{**} \in S$.
Now, from Theorem~\ref{thm:SA-JSQ_Process_level_convergence}, we have that the distribution of $\mf x^{**}$ must be invariant under the map $\mf u \mapsto \mf x(t,\mf u)$.
By the global stability result proved earlier
it follows that the only measure invariant
under the map $\mf u \mapsto \mf x(t,\mf u)$
is $\delta_{x^*}$.
Hence, we must have $\mf x^{**}= \mf x^*$.

\section{Unique Characterisation of Stationary Distribution from Given State}
\label{app:Unique_charac}
The main idea of the proof is to show that the transition rates of the chain $(V_{l_1(\mf x),1},V_{l_2(\mf x),2},\dots,V_{l_M(\mf x),M})$ defined in~\eqref{eqn:tran_rates_m_dim_markov} have a form similar to that of the transition rates of the Markov chain defined in~\eqref{eqn:k-dim-chain}. The proof then follows by the application of Lemma~\ref{lem:M-dim-Markov}.

Consider the chain $$\mf C=(V_{l_1(\mf x),1},V_{l_2(\mf x),2},\dots,V_{l_M(\mf x),M}).$$
We first reduce the dimension of the chain above and then rearrange the remaining components by performing the following steps sequentially:
\begin{enumerate}
    \item {\bf Step 1}: We first note from~\eqref{eqn:tran_rates_m_dim_markov} that the stationary rate of the transition $V_{l_j(\mf x),j} \to V_{l_j(\mf x),j}-1$ is zero (whereas the rate of the transition $V_{l_j(\mf x),j} \to V_{l_j(\mf x),j}+1$ is strictly positive) for any component $j \in [M]$ for which either (i) $l_j(\mf x) >l_i(\mf x)$ for some $i<j$ or (ii) $l_j(\mf x)-1 >l_i(\mf x)$ for some $i\neq j$. This is because for pairs $i,j \in [M]$ satisfying any of the above two conditions we have $\pi_{\mf{x}}(V_{l_j(\mf x),i}=\infty)=1$ and if $V_{l_j(\mf x),i}=\infty$, then from the definition of $R_{i,j}$ we have $\indic{\mf V \in R_{l_j(\mf x),j}}=0$. Therefore, for each $j \in [M]$ which satisfies one of the above two conditions we have $\pi_{\mf x}(V_{l_j(\mf x),j}=\infty)=1$. Hence, in the first step, we remove all these components from the chain $\mf C$.
 \item {\bf Step 2}: Next, we arrange the remaining components of the chain $\mf C$ in the increasing order of their corresponding minimum queue lengths (i.e., increasing order of $l_j(\mf {x})$). The components having the same minimum queue length (i.e., the same value of $l_j(\mf x)$) are then arranged in the decreasing order of their service rates. It can then be easily verified from the definition of $R_{i,j}$ that the remaining components have transition rates of the form given by~\eqref{eqn:k-dim-chain}.
    \end{enumerate}
Let $\mf Y =(Y_{m_1',p_1},Y_{m_2',p_2},\dots,Y_{m_H',p_H})$ denote the chain obtained after performing the steps mentioned above, where  for each $i\in[H]$, $Y_{m_i',p_i}$ is the leftover component of the chain $\mf C$ with minimum queue length $m_i'$ in pool $p_i$. Note that $m_i'\leq 
m_{i+1}'$ and $p_i\leq p_{i+1}$ for all $i\in [H]$. Moreover, the transition rate from $Y_{m_i',p_i} \to Y_{m_i',p_i} +1$ is given by $\nu_i=\gamma_{p_i}\mu_{p_i}(x_{m_i',p_i}-x_{m_i'+1,p_i})$ and the rate of transition from $Y_{m_i',p_i} \to Y_{m_i',p_i} -1$ is $\lambda \indic{0=Y_{m_1',p_1}=\dots=Y_{m_{i-1}',p_{i-1}}< Y_{m_i',p_i}}$. Now it remains to verify that $\pi_{\mf x}(V_{m_i',p_i}=\infty)=1$ or $0$ for each $i\in[H]$ using Lemma~\ref{lem:M-dim-Markov}.

Define $\rho_i=\frac{\nu_i}{\lambda}$ for each $i\in [H]$ and $\rho_{H+1}=0$. To show that $\pi_{\mf x}(V_{m_i',p_i}=\infty)=1$ or $0$ for each $i\in[H]$ it is suffices to show that for each $0\leq L\leq H$ if $\sum_{i=1}^L \rho_i<1$ and $\sum_{i=1}^{L+1} \rho_i\geq 1$ then we have
\begin{align*}
    \pi_{\mf x}(V_{m_1',p_1}&=\infty)=\dots=\pi_{\mf x}(V_{m_L',p_L}=\infty)=0,\\
    \pi_{\mf x}(V_{m_{L+1}',p_{L+1}}&=\infty)=\dots=\pi_{\mf x}(V_{m_H',p_H}=\infty)=1.
\end{align*}
Suppose $\sum_{i=1}^L \rho_i<1$ and $\sum_{i=1}^{L+1} \rho_i\geq 1$. Now we use contradiction to prove that $\pi_{\mf x}(Y_{m_i',p_i}=\infty)=0$ for all $i\in [L]$.
Assume $\pi_{\mathbf{x}}(Y_{m_i',p_i}=\infty)=\epsilon\in (0,1]$ for $i\in[L]$. Also, assume $\bar{\pi}_{\mf x}$ to be the unique stationary distribution of $(Y_{m_1',p_1},\dots,Y_{m_i',p_i})$ given that $Y_{m_r',p_r}$ is finite for all $r\in[i]$. Note that the stationary measure $\bar{\pi}_{\mf x}$ is same as the stationary measure $\pi$ given in Lemma~\ref{lem:M-dim-Markov}.  Therefore, we can write
\begin{align*}
\pi_{\mathbf{x}}(R_{m_i',p_i})&=(1-\epsilon)\bar{\pi}_{\mf x} (0=Y_{m_1',p_1}=\dots=Y_{m_{i-1}',p_{i-1}}< Y_{m_i',p_i})+\epsilon\\
&=(1-\epsilon)\frac{\nu_i}{\lambda} +\epsilon,    
\end{align*}
where we use $\bar{\pi}_{\mf x} (0=Y_{m_1',p_1}=\dots=Y_{m_{i-1}',p_{i-1}}< Y_{m_i',p_i})=\frac{\nu_i}{\lambda}$ from~\eqref{eqn:stat_prop_finite_M}. Now substituting the above in the differential form of~\eqref{eqn:SA-JSQ_fluid_process} we obtain
\eq{
\begin{aligned}
\frac{dx_{m_i',p_i}(t)}{dt}&= \frac{\lambda}{\gamma_{p_1}}\Big((1-\epsilon)\frac{\nu_i}{\lambda} +\epsilon\Big) -\mu_{p_i}(x_{m_1',p_1}(t)-x_{m_i'+1,p_i}(t))\\
&=\epsilon \Big(\frac{\lambda}{\gamma_{p_i}}-\mu_{p_1}(x_{m_i',p_i}(t)-x_{m_i'+1,p_i}(t))\Big)>0,
\end{aligned}
}
where the last step follows as $\frac{\nu_i}{\lambda}<1$. Now, since $x_{m_i',p_i}(t)=1$ as $m_i'$ is the minimum queue length in pool $p_i$, we must have $\frac{dx_{m_i',p_i}(t)}{dt}<0$ which leads to a contradiction. Therefore, we have $\pi_{\mathbf{x}}(Y_{m_i',p_i}=\infty)=0$ for all $i\in[L]$.
Now observe that $\sum_{i=1}^{L+1} \rho_i\geq 1$ implies $\rho_{L+1} \geq 1- \sum_{i=1}^{L}\rho_i$. Therefore, using Lemma~\ref{lem:M-dim-Markov}, we know that the chain $(Y_{m_1',p_1},\dots,Y_{m_{L+1}',p_{L+1}})$ is unstable. Moreover, again from Lemma~\ref{lem:M-dim-Markov}, $\sum_{i=1}^{L}\rho_i<1$ insures that the chain $(Y_{m_1',p_1},\dots,Y_{m_{L}',p_{L}})$ is stable. Hence, we have $\bar{\pi}_{\mathbf{x}}(Y_{m_{L+1}',p_{L+1}}\geq l)=1$ for all $l\geq0$. This shows that $\pi_{\mathbf{x}}(Y_{m_{L+1}',p_{L+1}}=\infty)=1$. Now observe that the rate of transition from $Y_{m_i',p_i} \to Y_{m_i',p_i}-1 $ is $\lambda \indic{0=Y_{m_1',p_1}=\dots=Y_{m_{i-1}',p_{i-1}}< Y_{m_i',p_i}}=0$ for each $i\in \cbrac{L+2,\dots,H}$. Hence, with only non-zero rate of transition from $Y_{m_i',p_i} \to Y_{m_i',p_i}+1$ we have $\pi_{\mathbf{x}}(Y_{m_{i}',p_{i}}=\infty)=1$ for each $i\in \cbrac{L+2,\dots,H}$.

\end{document}